\documentclass[12pt]{amsart}
\usepackage[centertags]{amsmath}
\usepackage{amstext,amssymb,amsopn,amsthm}
\usepackage{enumerate}
\usepackage[bookmarks,colorlinks]{hyperref}
\usepackage{color}
\usepackage[margin=1.5in]{geometry}
\usepackage[english]{babel}
\newcommand{\RR}{\mathbb{R}}
\newcommand{\NN}{\mathbb{N}}
\newcommand{\PP}{\mathbb{P}}
\newcommand{\EE}{\mathbb{E}}
\newcommand{\DD}{\mathbb{D}}
\newcommand{\HH}{\mathbb{H}}
\newcommand{\LL}{\mathbf{L}}
\newcommand{\s}{\mathbf{S}}
\newcommand{\Y}{\mathbf{Y}}
\newcommand{\p}{\mathcal{P}}
\newcommand{\h}{\mathcal{H}}
\newcommand{\f}{\mathcal{F}}
\newcommand{\M}{\mathcal{M}}
\newcommand{\C}{\mathcal{C}}
\newcommand{\A}{\mathcal{A}}
\newcommand{\T}{\mathcal{T}}
\newcommand{\x}{\overline{x}}
\newcommand{\y}{\overline{y}}

\newtheorem{theorem}{Theorem}[section]
\newtheorem{lemma}[theorem]{Lemma}
\newtheorem{corollary}[theorem]{Corollary}
\newtheorem{proposition}[theorem]{Proposition}

\theoremstyle{definition}
\newtheorem{definition}[theorem]{Definition}
\newtheorem{remark}[theorem]{Remark}

\begin{document}

\title[Boundary behavior of $\alpha$-harmonic functions]{Boundary behavior of $\alpha$-harmonic functions on the complement of the sphere and hyperplane}

\author{Tomasz Luks}
\address{
LAREMA, Laboratoire de Math\'ematiques\\Universit\'e d'Angers\\
2 boulevard Lavoisier\\
49045 Angers Cedex 01, France 
}
\email{luks@math.univ-angers.fr}
\thanks{This research was partially supported by Agence Nationale de la Recherche grant ANR-09-BLAN-0084-01 and by MNiSW grant N N201 373136.}

\keywords{$\alpha$-harmonic functions, fractional Laplacian, Hardy spaces, stable L\'evy process}
\subjclass[2010]{Primary 60J75; Secondary 60J50, 60J45, 42B30, 31B25.}

\begin{abstract}
We study $\alpha$-harmonic functions on the complement of the sphere and on the complement of the hyperplane 
in Euclidean spaces of dimension bigger than one, for $\alpha\in(1,2)$. We describe the corresponding  
Hardy spaces and prove the Fatou theorem for $\alpha$-harmonic functions. We also give explicit formulas 
for the Martin kernel of the complement of the sphere and for the harmonic measure, Green function and Martin kernel
of the complement of the hyperplane for the symmetric $\alpha$-stable L\'evy processes. Some extensions 
for the relativistic $\alpha$-stable processes are discussed.
\end{abstract}

\maketitle

\section{Introduction}\label{sec1}

Let $X_t$ be a symmetric $\alpha$-stable L\'evy process on $\RR^d$, $d\geq2$, with the index $\alpha\in(0,2)$ 
and the characteristic function
\begin{equation}\label{det29}
\EE^xe^{i\xi\cdot(X_t-x)}=e^{-t|\xi|^{\alpha}},\quad x,\xi\in\RR^d,\quad t\geq0.
\end{equation}
Here $\EE^x$ is the expectation for the process starting from $x$ and $\cdot$ denotes the standard inner product.
The study of the {\it $\alpha$-harmonic} functions, i.e., functions which are harmonic for $X_t$ (see Preliminaries for the definition), has 
been of interest in recent years, see \cite{Bog, CS1, Bog2, BogB, MS, SW, BaY, BD, MR, BJ, Ki, BogKK, BogZ, BBKRSV, C, MR1, BDL}. In some respects, 
the behavior of these functions contrasts sharply with that of the classical harmonic functions of the Laplacian (for which see, e.g., 
\cite{S1, S2, JK, ABR, A}). This is due to jumps of the process $X_t$. 
One of the most important properties distinguishing $X_t$ from diffusions is the fact that $X_t$ does not hit the boundary while 
leaving a sufficiently regular domain. Instead it jumps to the interior of the complement of the domain. 
Discontinuity of the trajectories at the exit time has a significant influence on the 
boundary properties of $\alpha$-harmonic functions. For example, it allows the existence of a positive $\alpha$-harmonic 
function on the unit ball in $\RR^d$ with the boundary limits identically equal to $\infty$, see \cite[Example 3.3., p. 59]{BBKRSV}. In the 
classical case of the Laplacian, when the underlying process is the Brownian motion, such examples do not 
exist because of the Fatou theorem and the so-called nontangential convergence of positive harmonic functions, see \cite{W,S1,S2,An,Wu1,JK,Do2,ABR,A}. 
In fact, Fatou-type theorems for $\alpha$-harmonic functions require an appropriate normalization, see \cite{BD, MR, Ki, MR1} and \cite{BaY}. 

In this paper we show that better analogues of the classical theory are obtained for sets which are complements of {\it smooth surfaces}, 
rather than for smooth domains. We proceed by considering two particular examples: the surface of the sphere $\s:=\left\{x\in\RR^d:|x|=1\right\}$ 
and of the hyperplane $\LL:=\left\{x=(x_1,...,x_d)\in\RR^d:x_d=0\right\}$. Here $d\geq2$ and $\alpha\in(1,2)$.

We should note that $\s$ and $\LL$ are non-polar sets for $X_t$ whenever $\alpha\in(1,2)$. The hitting probability 
of $\s$ for $X_t$ is given in \cite[Theorem 2.1]{P2}. Correspondingly, the last coordinate of $X_t=(X^1_t,...,X^d_t)$ is 
a one-dimensional symmetric $\alpha$-stable L\'evy motion, which 
hits zero almost surely for $\alpha\in(1,2)$ (i.e., $X^d_t$ is pointwise recurrent, see \cite{P2}).
We give the explicit formula for the hitting distribution of $\LL$ for $X_t$ in Section~\ref{sec4}, 
see (\ref{det32}) and Proposition~\ref{prop:poisker}. 
It turns out that the trajectories of $X_t$ are almost surely continuous on hitting $\s$ and $\LL$,
as in the case of the Brownian motion. We would like to recall the following fact. If $u$ is a nonnegative 
$\alpha$-harmonic function on an open set $D$ and $\tau_{D}$ is the first exit time 
from $D$ for $X_t$, then $M_t=u(X_{t\wedge\tau_{D}})$ is a positive supermartingale. By Doob's theory, 
$M_t$ has limits as $t\to\infty$. If $X_t$ is continuous at $t=\tau_D$, then 
we can translate the convergence of $M_t$ into the existence of the nontangential limits for $u$. Such approach has been proposed by
Doob in \cite{Do2} in order to prove the classical Fatou theorem and serves us as a motivation for studying $\DD:=\RR^d\setminus\s$ and 
$\HH:=\RR^d\setminus\LL$. The corresponding $\alpha$-harmonic versions of the Fatou theorem on $\DD$ and on $\HH$ are given 
as Theorem~\ref{th:Fatou} and Theorem~\ref{th:Fatou2}. We should recall that the result does not hold for nonnegative $\alpha$-harmonic 
functions on bounded domains, even for domains with very regular boundary (\cite[Example 3.3., p. 59]{BBKRSV}).

We are also interested in the Hardy spaces of $\alpha$-harmonic functions on $\DD$ and on $\HH$. 
The topic was intensively studied for classical harmonic functions, mainly on the ball 
and on the half-space, but also on bounded smooth and Lipschitz domains, see 
\cite{D, S1, S2, JK, K, ABR} (wider classes of diffusion operators were considered in \cite{W, L}). 
The classical Hardy spaces on the unit ball in $\RR^d$ are defined by the condition 
\[
\sup_{0\leq r<1}\int_{\s}|u(rx)|^p\sigma(dx)<\infty,
\]
where $\sigma$ is the surface measure. We consider analogous (analytic) definitions of Hardy spaces of $\alpha$-harmonic functions 
for $\DD$ and $\HH$. The analytic Hardy 
spaces are denoted by $h^p_{\alpha}(\DD)$ and $h^p_{\alpha}(\HH)$, respectively (see Definitions~\ref{def:hpan1} and~\ref{def:hpan2}), 
and they are characterized in Theorems~\ref{th:hpnorm},~\ref{th:spherehp} and~\ref{th:hyperplhp}. We also discuss the probabilistic 
version of Hardy spaces, denoted by $\h^p_{\alpha}(\DD)$ and $\h^p_{\alpha}(\HH)$, respectively (Definitions~\ref{def:hpprob1} 
and~\ref{def:hpprob2}). We refer the reader to \cite{BDL}, where the spaces were first described by the so-called Hardy-Stein 
identity for arbitrary open sets in $\RR^d$. We show that the analytic and the probabilistic definitions are equivalent on $\DD$ for 
all $p\in[1,\infty)$, see Lemma~\ref{lem:harmaj} and Theorem~\ref{th:hpprob1}, while the space $\h^p_{\alpha}(\HH)$ is essentially 
bigger than $h^p_{\alpha}(\HH)$, see Theorem~\ref{th:widerhp}. 

In the paper we also prove the following new explicit formulas: 
for the Martin kernel of $\DD$ for $X_t$ (Proposition~\ref{prop:martin}), for the hitting distribution of $\LL$ for $X_t$ 
(Proposition~\ref{prop:poisker}; the formula has been proved in \cite{I} only for $d=2$), and for the Green function 
and the Martin kernel of $\HH$ for $X_t$ (Proposition~\ref{prop:greenmartin}). The existence of such formulas was one of 
the main motivations for this work. Another motivation is to encourage the study of the boundary value problems for the 
fractional Laplacian $\Delta^{\alpha/2}$. Our study shows that the boundary conditions of Dirichlet type are of substantial 
interest when defined on {\it smooth surfaces} for $\alpha\in(1,2)$.

When $\alpha\in(0,1]$, $\s$ and $\LL$ are polar for $X_t$ 
(the hitting probability is zero), so the above mentioned probabilistic motivation disappears. However, it is possible to 
consider the conditional Hardy spaces in the spirit of \cite{MR1} and \cite{BDL}.

In the last part of the paper we show that similar problems can also be studied for other types of processes with discontinuous 
trajectories. As an example we give the so-called relativistic $\alpha$-stable process $X^m_t$, where $m>0$ is a parameter 
(see Section~\ref{sec5} for details). We prove that the sphere and the hyperplane are non-polar for $X^m_t$ if and only if 
$\alpha\in(1,2)$ (Proposition~\ref{prop:sphererelat} and Remark~\ref{rem:hyperrelat}). We also give the explicit formula for 
the so-called $\lambda$-harmonic measure of $\HH$ for $X^m_t$ in the particular case $\lambda=m$, $\alpha\in(1,2)$ 
(Proposition~\ref{prop:hyperplrelat}).

The paper is organized as follows. In Section~\ref{sec2} we give basic definitions and facts concerning the 
potential theory of the symmetric $\alpha$-stable processes and we recall the main results of \cite{P2}: the hitting probability of $\s$, 
the hitting distribution of $\s$ and the Green function of $\DD$ for $X_t$. In Section~\ref{sec3} we prove the explicit formula for 
the Martin kernel of $\DD$ for $X_t$ and we characterize the structure of the $\alpha$-harmonic Hardy spaces on $\DD$. We also prove 
the Fatou theorem for $\alpha$-harmonic functions. Analogous problems for $\HH$ are discussed in Section~\ref{sec4}. 
In Section~\ref{sec5} we analyze the case of the relativistic stable processes.

\section{Preliminaries}\label{sec2}

Throughout the paper, $(X_t,\PP^x)$ will denote a symmetric $\alpha$-stable L\'evy process on $\RR^d$ starting from $x$ and given 
by the characteristic function (\ref{det29}). We will consider only the index $\alpha\in(1,2)$ and the dimension $d\geq2$ 
unless stated otherwise. For a set $B\subset\RR^d$ let $\partial B$ denote the boundary of $B$, and let 
$\tau_B=\inf\left\{t>0:X_t\notin B\right\}$, $T_B=\tau_{B^c}$ be the first exit time and the first entry time for $B$, respectively. 
In the paper we use a convention that constants denoted by small letters may differ in each lemma, while constants denoted by 
capital letters do not change. The notation $c=c(a,b)$ means that the constant $c$ depends only on $a$ and $b$.

Let $D\subset\RR^d$ be open. A Borel function $u\colon\RR^d\to\RR$ is called $\alpha$-harmonic on $D$ if for every bounded open 
$U\subset D$ with $\overline{U}\subset D$ (denoted $U\subset\subset D$) we have
\begin{equation}\label{det30}
u(x)=\EE^x u(X_{\tau_U}),\quad x\in U.
\end{equation}
If $u\equiv0$ on $D^c$, then $u$ is called {\it singular} $\alpha$-harmonic on $D$. If (\ref{det30}) holds with $U=D$, then $u$ 
is called {\it regular} $\alpha$-harmonic on $D$. By the strong Markov property, every regular $\alpha$-harmonic function on $D$ 
is also $\alpha$-harmonic on $D$.

Equivalently, a Borel function $u$ on $\RR^d$ is $\alpha$-harmonic on $D$ if it is continuous on $D$, 
$\int_{\RR^d} |u(y)|(1+|y|)^{-d-\alpha}dy<\infty$, and
\begin{equation}\label{det27}
 \Delta^{\alpha/2}u(x)=\lim_{\varepsilon\to0}\A_{d,-\alpha}\int_{|x-y|>\varepsilon}\frac{u(y)-u(x)}{|x-y|^{d+\alpha}}dy=0,\quad x\in D,
\end{equation}
see \cite[Theorem 3.9]{BogB}. Here $\A_{d,\gamma}=\Gamma((d-\gamma)/2)/(2^{\gamma}\pi^{d/2}|\Gamma(\gamma/2)|)$ for 
$-2<\gamma<2$. The operator $\Delta^{\alpha/2}$ is called the fractional Laplacian and it is the infinitesimal generator of the process $X_t$.

For $x\in D$, the {\it $\alpha$-harmonic measure} for $D$ is the measure $\omega^x(\cdot,D)$ on $D^c$ given by
\begin{equation}\label{det31}
\omega^x(A,D):=\PP^x(X_{\tau_D}\in A;\tau_D<\infty),\quad A\subseteq D^c,
\end{equation}
i.e., $\omega^x(\cdot,D)$ is the exit distribution from $D$ for $X_t$ starting from $x$. The Green function of $D$ for $X_t$ is defined by
\[
 G_D(x,y)=\A_{d,\alpha}\left(|x-y|^{\alpha-d}-\int_{D^c}|y-z|^{\alpha-d}\omega^x(dz,D)\right),\quad x,y\in\RR^d.
\] 
The Poisson kernel of $D$ for $X_t$, $\overline{D}\neq\RR^d$, is given by the formula
\begin{equation}\label{det24}
 P_D(x,y)=\A_{d,-\alpha}\int_D\frac{G_D(x,z)}{|z-y|^{d+\alpha}}dz,\quad x\in D, y\in D^c.
\end{equation}
We have
\[
 \omega^x(A,D)=\int_AP_D(x,y)dy,\quad x\in D, A\subseteq(\overline{D})^c,
\]
i.e., $P_D(x,y)$ is the density of $\omega^x(\cdot,D)$ with respect to the Lebesgue measure on $(\overline{D})^c$ 
and corresponds to the effect of leaving $D$ by a jump (see \cite{Wu2} for details). 
If $D=B(a,r)$ is a ball of center $a\in\RR^d$ and radius $r>0$, then we have
\begin{equation}\label{det23}
 P_{B(a,r)}(x,y)=\C_1\left(\frac{r^2-|x-a|^2}{|y-a|^2-r^2}\right)^{\alpha/2}\frac{1}{|x-y|^d},
\end{equation}
where $\C_1=\Gamma(d/2)\pi^{-1-d/2}\sin(\pi\alpha/2)$. For $D=\DD$ or $D=\HH$ we have $\overline{D}=\RR^d$ and 
the right hand side of (\ref{det24}) is equal to $\infty$ for all $x\in D$ and $y\in D^c$ (see the proof of 
Propositions~\ref{prop:Martin1}).

By \cite[Theorem 2.1]{P2}, the hitting probability $\Phi(x):=\PP^x(T_{\s}<\infty)$ is given by
\begin{equation}\label{det3}
\Phi(x)=\begin{cases}
		\C_2\left||x|^2-1\right|^{(\alpha/2)-1}|x|^{1-d/2}P^{1-d/2}_{-\alpha/2}\left(\frac{|x|^2+1}{\left||x|^2-1\right|}\right), & \ x\in\DD,x\neq0,\\
		\C_2/\Gamma(d/2), & \ x=0,
		\end{cases}
\end{equation}
where $\C_2=\pi^{1/2}2^{2-\alpha}\Gamma((\alpha+d)/2-1)/\Gamma((\alpha-1)/2)$ and $P^{\nu}_{\mu}$ is the usual 
Legendre function of the first kind. Since $\Phi$ is radial, for $r\geq0$, $r\neq1$ we let $\phi(r):=\Phi(rz)$, $z\in\s$. 
By \cite[Vol. I, formula 20, p. 164]{E} we 
have $\phi(r)\to1$ when $r\to1$ and from \cite[Vol. I, formula 3, p. 163]{E} it follows that $\phi(r)\to0$ as $r\to\infty$.

Let $\sigma$ be the normalized surface measure on $\s$. By \cite[Theorem 3.1]{P2}, the hitting distribution of $\s$ for $X_t$ 
(or the $\alpha$-harmonic measure for $\DD$) has a density with respect to $\sigma$ given by
\begin{equation}\label{det1}
\p_{\DD}(x,y)=\frac{\C_2}{\Gamma(d/2)}\frac{\left||x|^2-1\right|^{\alpha-1}}{\left|x-y\right|^{d+\alpha-2}},\quad x\in\DD, y\in\s.
\end{equation}
In particular, $\Phi(x)=\int_{\s}\p_{\DD}(x,y)\sigma(dy)$. We will call $\p_{\DD}(x,y)$ the Poisson kernel 
of $\DD$ for $X_t$. By the symmetry we have
\begin{equation}\label{det2}
\p_{\DD}(ry,z)=\p_{\DD}(rz,y),\quad y,z\in\s ,r\geq0,r\neq1.
\end{equation}
By \cite[Theorem 4.1]{P2}, the Green function of ${\DD}$ for $X_t$ is given by
\begin{equation}\label{det5}
G_{\DD}(x,y)=\frac{\A_{d,\alpha}}{|x-y|^{d-\alpha}}\left[1-\Phi\left(\frac{y}{|y-x|}\left|x-y/|y|^2\right|\right)\right],\quad x,y\in\DD.
\end{equation}

For an open set $D\subset\RR^d$ and a nonnegative function $f$ on $\RR^d$ we define
\begin{equation}\label{det22}
\f_D[f](x)=\sup_{U\subset\subset D}\EE^xf(X_{\tau_U}).
\end{equation}
The function $\f$ will play an important role in our study. Basically, we will consider only $\f_{\DD}$ and $\f_{\HH}$, 
but the following result remains true for arbitrary open subsets of $\RR^d$.

\begin{lemma}\label{lem:majorant1}
Let $D\subset\RR^d$ be open and let $u$ be $\alpha$-harmonic on $D$. Suppose that for some $x\in D$ and some $p\in[1,\infty)$ we have 
$\f_D[|u|^p](x)<\infty$. Then $\f_D[|u|^p]$ is the minimal $\alpha$-harmonic majorant of $|u|^p$ on $D$. 
\end{lemma}

\begin{proof}
Let $U_n$ be a sequence of open bounded sets such that $\overline{U}_n\subset D$, $U_n\subseteq U_{n+1}$ for every $n$ and $\bigcup_nU_n= D$. Set 
$\tau_n:=\tau_{U_n}$. Then $\tau_n\leq\tau_{n+1}$, so by the strong Markov property and the Jensen inequality we have
\[
\EE^x|u(X_{\tau_n})|^p\leq\EE^x[\EE^{X_{\tau_{n}}}|u(X_{\tau_{n+1}})|^p]=\EE^x|u(X_{\tau_{n+1}})|^p.
\]
Hence $\f_D[|u|^p](x)=\lim_n\EE^x|u(X_{\tau_n})|^p$. The monotone convergence theorem and the Harnack inequality (see \cite[Theorem 1]{BogKK}) 
imply that either $\f_D[|u|^p]$ is $\alpha$-harmonic on $\DD$ or $\mathcal{F}_{\DD}[|u|^p]\equiv\infty$. Therefore, if 
$\f_D[|u|^p](x)<\infty$ for some $x\in D$ and some $p\in[1,\infty)$ then $\f_D[|u|^p]$ is $\alpha$-harmonic 
and nonnegative on $D$. By Jensen's inequality we have $|u|^p\leq\EE^x|u(X_{\tau_n})|^p\leq\f_D[|u|^p]$. Furthermore, if $h$ is nonnegative 
and $\alpha$-harmonic on $D$ such that $|u|^p\leq h$, then $\EE^x|u(X_{\tau_n})|^p\leq\EE^xh(X_{\tau_n})=h(x)$ for every $n$. Hence 
$\f_D[|u|^p]\leq h$, as desired.
\end{proof}

\section{$\alpha$-harmonic functions on the complement of the sphere}\label{sec3}

In this section we will characterize the behavior of $\alpha$-harmonic functions on $\DD=\RR^d\setminus\s$. 
We denote by $C(\s)$ the space of continuous functions on $\s$ and by $\M(\s)$ the space of finite 
signed Borel measures on $\s$. For $\mu\in\M(\s)$ let $\|\mu\|$ denote the total variation norm of $\mu$. 
For $1\leq p<\infty$ and a Borel function $f$ on $\s$ let 
\[
\|f\|_p:=\left(\int_{\s}|f(x)|^p\sigma(dx)\right)^{1/p},
\]
and let $\|f\|_{\infty}$ denote the essential supremum norm on $\s$ with respect to $\sigma$. 

For $f\in L^1(\s,\sigma)$ and $\mu\in\M(\s)$, we define the {\it Poisson integrals} of $f$ and $\mu$ on $\DD$ as
\[
 \p_{\DD}[\mu](x)=\int_{\s}\p_{\DD}(x,y)\mu(dy),\quad \p_{\DD}[f](x)=\int_{\s}\p_{\DD}(x,y)f(y)\sigma(dy).
\]
For every $f\in L^1(\s,\sigma)$, the function $u$ given by 
\[
u(x)=\begin{cases}
		\p_{\DD}[f](x), & \ x\in\DD,\\
		f(x), & \ x\in\s,
		\end{cases}
\]
is regular $\alpha$-harmonic on $\DD$, and hence $\p_{\DD}[f]$ is $\alpha$-harmonic on $\DD$. For every $\mu\in\M(\s)$, $\p_{\DD}[\mu]$ is also
$\alpha$-harmonic on $\DD$, what can be shown by considering a sequence of continuous functions on $\s$ which converges to $\mu$ in weak$^*$ 
topology. See also Proposition~\ref{prop:Martin1}.

We define the Martin kernel of $\DD$ for $X_t$ as 
\begin{equation}\label{det6}
M_{\DD}(x,z)=\lim_{\DD\ni y\to z}\frac{G_{\DD}(x,y)}{G_{\DD}(0,y)},\quad x\in\DD,z\in\s\cup\left\{\infty\right\}.
\end{equation}
By \cite[Theorem 2]{BogKK}, the limit in (\ref{det6}) always exists. Directly from (\ref{det5}) we get
\begin{equation}\label{det9}
M_{\DD}(x,\infty)=\frac{1-\Phi(x)}{1-\Phi(0)}.
\end{equation}
For $z\in\s$ the formula is given in the next proposition, and it shows the relation between the Martin kernel and the $\alpha$-harmonic 
measure of $\DD$. We would like to remark that, probably, this observation may be generalized for wider class of open sets with 
non-polar boundary. 

\begin{proposition}\label{prop:martin}
For all $x\in\DD$ and $z\in\s$ we have
\begin{equation}\label{det8}
M_{\DD}(x,z)=\frac{\p_{\DD}(x,z)}{\p_{\DD}(0,z)}=\frac{\left||x|^2-1\right|^{\alpha-1}}{\left|x-z\right|^{d+\alpha-2}}.
\end{equation}
\end{proposition}

\begin{proof}
We show the proposition only for $|x|<1$, for $|x|>1$ the proof is similar. Since the limit in (\ref{det6}) exists for every $z\in\s$, we may assume that $y=rz$, and take $r\to1^-$. Let
\[
w_r=\frac{rz}{|rz-x|}\left|x-z/r\right|.
\]
Then $|w_r|>1$, $w_r\to z$ as $r\to1^-$ and by (\ref{det5}) we have
\[
M_{\DD}(x,z)=\frac{1}{|x-z|^{d-\alpha}}\cdot\lim_{r\to1^-}\frac{1-\Phi\left(w_r\right)}{1-\Phi\left(z/r\right)}.
\]
By \cite[Vol. I, formula (9), p. 123 and formula (18), p. 125]{E} we have the following expansion of the Legendre function
\[
P^{1-d/2}_{-\alpha/2}(t)=f_1(d,\alpha,t)F\left(1-\frac{\alpha}{2},\frac{d-\alpha}{2};2-\alpha;\frac{2}{1+t}\right) 
\]
\[
+ f_2(d,\alpha,t)F\left(\frac{\alpha}{2},\frac{d+\alpha}{2}-1;\alpha;\frac{2}{1+t}\right),\quad t>1,
\]
where $F$ is the hypergeometric function given by
\[
F(a,b;c;s)=\sum^{\infty}_{n=0}\frac{(a)_n(b)_n}{(c)_nn!}s^n,\quad c\neq0,-1,-2,...,
\]
$(\cdot)_n:=\Gamma(\cdot+n)/\Gamma(\cdot)$, and
\[
f_1(d,\alpha,t)=\frac{2^{1-\alpha/2}\Gamma(\alpha-1)}{\Gamma\left(\frac{\alpha}{2}\right)\Gamma\left(\frac{\alpha+d}{2}-1\right)}\cdot(t+1)^{\alpha/2-d/4-1/2}(t-1)^{d/4-1/2},
\]
\[
f_2(d,\alpha,t)=\frac{2^{\alpha/2}\Gamma(1-\alpha)}{\Gamma\left(1-\frac{\alpha}{2}\right)\Gamma\left(\frac{d-\alpha}{2}\right)}
\cdot(t+1)^{1/2-d/4-\alpha/2}(t-1)^{d/4-1/2}.
\]
For $v\in D$, $|v|>1$ denote
\[
I(d,\alpha,v)=\C_2\left(|v|^2-1\right)^{\alpha/2-1}|v|^{1-d/2}.
\]
Then we have
\[
\Phi(v)=I(d,\alpha,v)P^{1-d/2}_{-\alpha/2}\left(\frac{|v|^2+1}{|v|^2-1}\right)
\]
\[
=I(d,\alpha,v)\left[f_1\left(d,\alpha,\frac{|v|^2+1}{|v|^2-1}\right)+f_2\left(d,\alpha,\frac{|v|^2+1}{|v|^2-1}\right)
+G(d,\alpha,v)\right],
\]
where
\[
G(d,\alpha,v)=f_1\left(d,\alpha,\frac{|v|^2+1}{|v|^2-1}\right)\sum^{\infty}_{n=1}A_n(d,\alpha)\left(\frac{|v|^2-1}{|v|^2}\right)^n
\]
\[
+f_2\left(d,\alpha,\frac{|v|^2+1}{|v|^2-1}\right)\sum^{\infty}_{n=1}B_n(d,\alpha)\left(\frac{|v|^2-1}{|v|^2}\right)^n,
\]
and
\[ 
A_n(d,\alpha)=\frac{\left(1-\frac{\alpha}{2}\right)_n\left(\frac{d-\alpha}{2}\right)_n}{(2-\alpha)_nn!},\quad
B_n(d,\alpha)=\frac{\left(\frac{\alpha}{2}\right)_n\left(\frac{d+\alpha}{2}-1\right)_n}{(\alpha)_nn!}.
\]
We have
\[
f_1\left(d,\alpha,\frac{|v|^2+1}{|v|^2-1}\right)=\frac{\Gamma(\alpha-1)}{\Gamma\left(\frac{\alpha}{2}\right)\Gamma\left(\frac{\alpha+d}{2}-1\right)}\cdot|v|^{\alpha-1-d/2}\left(|v|^2-1\right)^{1-\alpha/2},
\]
\[
f_2\left(d,\alpha,\frac{|v|^2+1}{|v|^2-1}\right)=\frac{\Gamma(1-\alpha)}{\Gamma\left(1-\frac{\alpha}{2}\right)\Gamma\left(\frac{d-\alpha}{2}\right)}\cdot|v|^{1-d/2-\alpha}\left(|v|^2-1\right)^{\alpha/2}.
\]
Furthermore, using the duplication formula,
\[
I(d,\alpha,v)f_1\left(d,\alpha,\frac{|v|^2+1}{|v|^2-1}\right)=|v|^{\alpha-d}.
\]
We write
\[
I(d,\alpha,v)f_2\left(d,\alpha,\frac{|v|^2+1}{|v|^2-1}\right)=c\left(|v|^2-1\right)^{\alpha-1}|v|^{2-d-\alpha},
\]
where $c=c(d,\alpha)=\C_22^{\alpha/2}\Gamma(1-\alpha)/[\Gamma\left(1-\alpha/2\right)\Gamma((d-\alpha)/2)]$. Hence we obtain
\[
 1-\Phi(v)=1-|v|^{\alpha-d}-c\left(|v|^2-1\right)^{\alpha-1}|v|^{2-d-\alpha}-I(d,\alpha,v)G(d,\alpha,v)
\]
\[
=\left(|v|^2-1\right)^{\alpha-1}\left[\frac{1-|v|^{\alpha-d}}{\left(|v|^2-1\right)^{\alpha-1}}
-c|v|^{2-d-\alpha}-\frac{I(d,\alpha,v)G(d,\alpha,v)}{\left(|v|^2-1\right)^{\alpha-1}}\right].
\]
Since $\alpha<2$,
\[
\lim_{|v|\to1^+}\frac{1-|v|^{\alpha-d}}{\left(|v|^2-1\right)^{\alpha-1}}=0.
\]
Furthermore,
\[
\frac{I(d,\alpha,v)G(d,\alpha,v)}{\left(|v|^2-1\right)^{\alpha-1}}=|v|^{\alpha-d-2}\left(|v|^2-1\right)^{2-\alpha}\sum^{\infty}_{n=1}A_n(d,\alpha)\left(\frac{|v|^2-1}{|v|^2}\right)^{n-1}
\]
\[
+c|v|^{2-d-\alpha}\sum^{\infty}_{n=1}B_n(d,\alpha)\left(\frac{|v|^2-1}{|v|^2}\right)^n,
\]
and thus
\[
\lim_{|v|\to1^+}\frac{I(d,\alpha,v)G(d,\alpha,v)}{\left(|v|^2-1\right)^{\alpha-1}}=0.
\]
We have
\[
|w_r|^2-1=\left|\frac{rz}{|rz-x|}\left|x-z/r\right|\right|^2-1=\frac{(1-r^2)\left(1-|x|^2\right)}{|rz-x|^2}.
\]
Therefore
\[
\lim_{r\to1^-}\frac{1-\Phi\left(w_r\right)}{1-\Phi\left(z/r\right)}=\lim_{r\to1^-}\left(\frac{|w_r|^2-1}{|z/r|^2-1}\right)^{\alpha-1}
=\lim_{r\to1^-}\left[\frac{(1-r^2)\left(1-|x|^2\right)}{|rz-x|^2(1/r^2-1)}\right]^{\alpha-1}
\]
\[
=\frac{\left(1-|x|^2\right)^{\alpha-1}}{|z-x|^{2\alpha-2}}\lim_{r\to1^-}\left(\frac{1-r^2}{1/r^2-1}\right)^{\alpha-1}
=\frac{\left(1-|x|^2\right)^{\alpha-1}}{|z-x|^{2\alpha-2}}.
\]
\end{proof} 

\noindent Proposition~\ref{prop:martin} shows that the Martin kernel and the Poisson kernel of $\DD$ are the same objects up to a multiplicative 
constants. This property together with the results of \cite{BogKK} give us the following representation for nonnegative 
$\alpha$-harmonic functions on $\DD$.

\begin{proposition}\label{prop:Martin1}
For every nonnegative measure $\mu\in\M(\s)$ and every constant $c\geq0$ the function $u$ given by
\begin{equation}\label{det16}
u(x)=\int_{\s}\p_{\DD}(x,y)\mu(dy)+c(1-\Phi(x)),\quad x\in\DD
\end{equation}
is $\alpha$-harmonic on $\DD$. Conversely, if $u$ is nonnegative and $\alpha$-harmonic on $\DD$ then there exists a unique nonnegative measure 
$\mu\in\M(\s)$ and a unique constant $c\geq0$ satisfying (\ref{det16}).
\end{proposition}

\begin{proof}
Since $\left(\overline{\DD}\right)^c=\emptyset$ and $\partial\DD=\s$ is of the Lebesgue measure 0, every $\alpha$-harmonic function on $\DD$ 
can be considered as a singular $\alpha$-harmonic function on $\DD$. The proposition is then a consequence of \cite[Lemma 14]{BogKK}, 
(\ref{det1}), (\ref{det9}) and (\ref{det8}) since all points of $\s\cup\left\{\infty\right\}$ are accessible from $\DD$ 
(see \cite[p.347]{BogKK} for the definition of accessibility). For $z\in\s$ this property can be deduced from \cite[(75) and 
(76)]{BogKK} and the fact, that $G_{B_z}(x,y)\leq G_{\DD}(x,y)$, where $B_z$ is a ball contained in $\DD$ and tangent to $\s$ at $z$. 
By (\ref{det24}) and (\ref{det23}) we then have 
\[
 \int_D\frac{G_{\DD}(0,y)}{|y-z|^{d+\alpha}}dy\geq P_{B_z}(0,z)=\infty.
\] 
On the other hand, from (\ref{det5}) it follows that $G_{\DD}(0,y)\approx\A_{d,\alpha}|y|^{\alpha-d}(1-\Phi(0))$ as $y\to\infty$. Hence $\int_{\RR^d}G_{\DD}(0,y)dy=\infty$, which gives the accessibility of the point at infinity.
\end{proof}

\noindent Proposition~\ref{prop:Martin1} implies that for every pair $(\mu,c)\in\M(\s)\times\RR$ the function $u(x)=\p_{\DD}[\mu](x)+c(1-\Phi(x))$ is 
$\alpha$-harmonic on $\DD$. This is a consequence of the Hahn decomposition $\mu=\mu^+-\mu^-$.

For $r>0$, $r\neq1$ and a function $u$ on $\DD$ we define the function $u_r$ on $\s$ by 
\begin{equation}\label{det33}
u_r(x):=u(rx),\quad x\in\s.
\end{equation}

\begin{lemma}\label{lem:poissinteg}
Poisson integrals on $\DD$ have the following properties:
\begin{enumerate}[\upshape (i)]
 \item If $\mu\in\M(\s)$, then $\|\p_{\DD}[\mu]_r\|_1\leq\|\mu\|$ for every $r$.
 \item If $1\leq p\leq\infty$ and $f\in L^p(\s,\sigma)$, then $\|\p_{\DD}[f]_r\|_p\leq\|f\|_p$ for every $r$.
 \item If $f\in C(\s)$, then $\|\p_{\DD}[f]_r-f\|_{\infty}\to0$ as $r\to1$.
 \item If $1\leq p<\infty$ and $f\in L^p(\s,\sigma)$, then $\|\p_{\DD}[f]_r-f\|_p\to0$ as $r\to1$.
 \item If $\mu\in\M(\s)$, then $\p_{\DD}[\mu]_r\to\mu$ weak$^*$ in $\M(\s)$ as $r\to1$.
 \item If $f\in L^{\infty}(\s,\sigma)$, then $\p_{\DD}[f]_r\to f$ weak$^*$ in $L^{\infty}(\s,\sigma)$ as $r\to1$.
\end{enumerate}
\end{lemma}

\begin{proof}
We start with the property (iii). Let $f\in C(\s)$. For $x\in\s$ and $r>0$, $r\neq1$ we have
\[
 |\p_{\DD}[f]_r(x)-f(x)|=\left|\int_{\s}\p_{\DD}(rx,y)f(y)\sigma(dy)-f(x)\right|
\]
\[
 \leq\int_{\s}\p_{\DD}(rx,y)|f(y)-f(x)|\sigma(dy)+(1-\phi(r))|f(x)|
\]
Let $\varepsilon>0$. Since $\s$ is compact, there exists $\delta>0$ independent of $x$ such that $|f(x)-f(y)|<\varepsilon$ if 
$|x-y|<\delta$. Hence the last term is less than
\[
 \varepsilon+\int_{|x-y|>\delta}\p_{\DD}(rx,y)|f(y)-f(x)|\sigma(dy)+(1-\phi(r))|f(x)|
\]
\[
 \leq\varepsilon+2\|f\|_{\infty}\int_{|x-y|>\delta}\p_{\DD}(x,y)\sigma(dy)+(1-\phi(r))\|f\|_{\infty}
\]
\[
\leq\varepsilon+\|f\|_{\infty}\left[2\C_1\Gamma(d/2)^{-1}|r^2-1|^{\alpha-1}\delta^{2-d-\alpha}+1-\phi(r)\right].
\]
The statement now follows from the fact that $\phi(r)\to1$ as $r\to1$ (see Preliminaries) and $\alpha\in(1,2)$. We prove the other statements of the lemma 
using the symmetry property (\ref{det2}) and the Jensen inequality, in a similar way as in the proofs of analogous properties in the classical case in 
\cite[Theorem 6.4, 6.7 and 6.9]{ABR}.
\end{proof}

\begin{corollary}\label{cor:Martinunique1}
Suppose $u$ is $\alpha$ harmonic on $\DD$. Then $u(x)=\p_{\DD}[\mu](x)+c(1-\Phi(x))$ for some pair $(\mu,c)\in\M(\s)\times\RR$ if and only 
if there exists a nonnegative $\alpha$-harmonic function $v$ on $\DD$ such that $|u|\leq v$. Furthermore, $\mu$ and $c$ are unique.
\end{corollary}

\begin{proof}
In view of Proposition~\ref{prop:Martin1}, we obtain the first part of the corollary in the same way as in the proof of 
\cite[Lemma 1]{MR1}. The uniqueness of the representation follows from Lemma~\ref{lem:poissinteg} (v).
\end{proof}

\begin{definition}\label{def:hpan1} For $p\in [1,\infty]$ we define the {\it Hardy space} $h^p_{\alpha}(\DD)$ as the family of functions $u$ 
$\alpha$-harmonic on $\DD$ such that
\[
\left\|u\right\|_{h^p}:=\sup_{r\in\RR_+\setminus\left\{1\right\}}\|u_r\|_p<\infty.
\]
\end{definition}

\noindent Note that since $\sigma(\s)<\infty$, for $1\leq p<q\leq\infty$ we have $h^q_{\alpha}(\DD)\subset h^p_{\alpha}(\DD)$, and 
\[
 \left\|u\right\|_{h^{\infty}}=\sup_{x\in\DD}|u(x)|.
\]
We will describe the spaces $h^p_{\alpha}(\DD)$ in terms of the Poisson integrals. The first part is the following.

\begin{theorem}\label{th:hpnorm}
Let $u$ be $\alpha$-harmonic on $\DD$.
\begin{enumerate}
	\item[\it 1.] If $u(x)=\p_{\DD}[\mu](x)+c(1-\Phi(x))$ for some pair $(\mu,c)\in\M(\s)\times\RR$, then 
$u\in h^1_{\alpha}(\DD)$ and $\|u\|_{h^1}=\|\mu\|\vee|c|$.
	\item[\it 2.] Let $p\in(1,\infty]$. If $u(x)=\p_{\DD}[f](x)+c(1-\Phi(x))$ for some pair $(f,c)\in L^p(\s,\sigma)\times\RR$, then 
$u\in h^p_{\alpha}(\DD)$ and $\|u\|_{h^p}=\|f\|_p\vee|c|$.
\end{enumerate}
\end{theorem}

\begin{proof}
Let $u(x)=\p_{\DD}[\mu](x)+c(1-\Phi(x))$ for some $(\mu,c)\in\M(\s)\times\RR$. By Lemma~\ref{lem:poissinteg} (i)
we have $\left\|\p_{\DD}[\mu]\right\|_{h^1}\leq\|\mu\|$. Since $\Phi$ is bounded, we have $u\in h^1_{\alpha}(\DD)$. 
From Lemma~\ref{lem:poissinteg} (v) it follows that $\|\mu\|\leq \liminf_{r\to1}\|\p_{\DD}[\mu]_r\|_1$, so
\[
\|\mu\|=\lim_{r\to1}\|\p_{\DD}[\mu]_r\|_1=\|\p_{\DD}[\mu]\|_{h^1}.
\]
We have $\|\p_{\DD}[\mu]_r\|_1\to0$ when $r\to\infty$. Furthermore, $\phi(r)\to1$ as $r\to1$ and $\phi(r)\to0$ 
as $r\to\infty$. Since 
\[
\|u_r\|_1\geq\left|\|\p_{\DD}[\mu]_r\|_1-|c|(1-\phi(r))\right|,
\]
we have $\|u\|_{h^1}\geq\|\mu\|\vee|c|$. On the other hand, by (\ref{det2}) we have
\[
\|u_r\|_1\leq\int_{\s}\int_{\s}\p_{\DD}(rx,y)|\mu|(dy)\sigma(dx)+|c|(1-\phi(r))
\]
\[
=\int_{\s}\int_{\s}\p_{\DD}(ry,x)\sigma(dx)|\mu|(dy)+|c|(1-\phi(r))
\]
\[
=\|\mu\|\phi(r)+|c|(1-\phi(r))\leq\|\mu\|\vee|c|.
\]
This gives the first part. To prove the second part choose $p\in(1,\infty]$ and suppose that $u(x)=\p_{\DD}[f](x)+c(1-\Phi(x))$ 
for some $(f,c)\in L^p(\s,\sigma)\times\RR$. Then by Lemma~\ref{lem:poissinteg} (ii) we have 
$\|\p_{\DD}[f]\|_{h^p}\leq\|f\|_p$, and since $\Phi$ is bounded 
we have $u\in h^p_{\alpha}(\DD)$. Using Lemma~\ref{lem:poissinteg} (iv) and (vi) we obtain
\[
\|\p_{\DD}[f]\|_{h^p}=\lim_{r\to1}\|\p_{\DD}[f]_r\|_p=\|f\|_p.
\]
Since 
\[
\|u_r\|_p\geq\left|\|\p_{\DD}[f]_r\|_p-|c|(1-\phi(r))\right|,
\]
as in the case $p=1$ we conclude that $\|u\|_{h^p}\geq\|f\|_p\vee|c|$. On the other hand, for $p<\infty$ we have
\[
\|u_r\|_p\leq\left(\int_{\s}\left|\int_{\s}\p_{\DD}(rx,y)f(y)\sigma(dy)\right|^p\sigma(dx)\right)^{\frac{1}{p}}+|c|(1-\phi(r)).
\]
By Jensen's inequality and (\ref{det2}), the last term above is less than
\[
\left(\phi(r)^{p-1}\int_{\s}\int_{\s}\p_{\DD}(ry,x)\sigma(dx)|f(y)|^p\sigma(dy)\right)^{\frac{1}{p}}+|c|(1-\phi(r))
\]
\[
=\phi(r)\|f\|_p+|c|(1-\phi(r))\leq\|f\|_p\vee|c|.
\]
For $p=\infty$ it suffices to estimate $f$ by $\|f\|_{\infty}$.
\end{proof}

\noindent Our aim now is to prove the inverse implication of Theorem~\ref{th:hpnorm}. In order to do this, we will first show that 
the spaces $h^p_{\alpha}(\DD)$ can be equivalently characterized by the function $\f_{\DD}$. For $n=2,3,...$ let 
\[
\DD_n=\left\{x:|x|<1-1/n\right\}\cup\left\{x:1+1/n<|x|<n\right\}.
\]
$\DD_n$ are bounded, $\overline{\DD}_n\subset\DD$, $\DD_n\subset\DD_{n+1}$ for every $n$ and $ \bigcup_n\DD_n=\DD$. 
The following property of the Poisson kernels $P_{\DD_n}(x,y)$ will be important in our approach.

\begin{lemma}\label{lem:rotation}
Fix $n\in\NN,n\geq2$. For all $r\in[0,1-1/n)\cup(1+1/n,n)$, $s\in(1-1/n,1+1/n)\cup(n,\infty)$ and $y\in\s$ we have
\begin{equation}\label{det38}
\int_{\s}P_{\DD_n}(rx,sy)\sigma(dx)=\int_{\s}P_{\DD_n}(ry,sx)\sigma(dx).
\end{equation}
Furthermore, the integrals are constant with respect to $y$. 
\end{lemma}

\begin{proof}
The identity (\ref{det38}) follows from the fact, that for $r,s$ as above and all $x,y\in\s$ we have 
$P_{\DD_n}(rx,sy)=P_{\DD_n}(ry,sx)$. This is a consequence of the rotation invariance and of the symmetry of the process $X_t$. Clearly, 
for any rotation $\mathcal{T}$ on $\RR^d$ and all $x,y\in\DD_n$ we have $G_{\DD_n}(x,y)=G_{\T\DD_n}(\T x,\T y)$. 
Obviously, $\T\DD_n=\DD_n$. Fix now $x,y\in\s$, $r\in[0,1-1/n)\cup(1+1/n,n)$ and $s\in(1-1/n,1+1/n)\cup(n,\infty)$. Let $\T_0$ be a rotation
on $\RR^d$ for which $\T_0x=y$. Then (\ref{det24}) and the rotation invariance imply that
\[
P_{\DD_n}(rx,sy)=P_{\DD_n}(r\T^{-1}_0x,sx).
\]
Furthermore, by the symmetry,
\[
P_{\DD_n}(r\T^{-1}_0x,sx)=P_{\DD_n}(r\T_0x,sx)=P_{\DD_n}(ry,sx).
\]
To prove the second part of the lemma, fix an arbitrary rotation $\T$. Then for $y\in\s$ and $r,s$ as before we have
\[
\int_{\s}P_{\DD_n}(rx,s\mathcal{T}y)\sigma(dx)=
\int_{\s}P_{\DD_n}(r\T^{-1}x,y)\sigma(dx)=\int_{\s}P_{\DD_n}(rz,sy)\sigma(dz).
\]
Clearly, $\sigma(\T(E))=\sigma(E)$ for every Borel set $E\subset\s$.
\end{proof}

\noindent As a consequence of Lemma~\ref{lem:rotation} we have the following equivalent characterization of the spaces 
$h^p_{\alpha}(\DD)$ for $p\in[1,\infty)$.

\begin{lemma}\label{lem:harmaj}
Let $u$ be $\alpha$-harmonic on $\DD$. Then $u\in h^p_{\alpha}(\DD)$ for a given $p\in[1,\infty)$ if and only if 
$\f_{\DD}[|u|^p](x)<\infty$ for some $x\in\DD$.  Furthermore, $\|u\|^p_{h^p}=\|\f_{\DD}[|u|^p]\|_{h^1}$.
\end{lemma}

\begin{proof}
Suppose first that $\f_{\DD}[|u|^p](x)<\infty$ for some $p\in[1,\infty)$ and some $x\in\DD$. Then by Lemma~\ref{lem:majorant1},
$\f_{\DD}[|u|^p]$ is nonnegative and $\alpha$-harmonic on $\DD$ and from Proposition~\ref{prop:Martin1} 
and Theorem~\ref{th:hpnorm} it follows that $\f_{\DD}[|u|^p]\in h^1_{\alpha}(\DD)$. Since 
$|u|^p\leq\f_{\DD}[|u|^p]$, we have $u\in h^p_{\alpha}(\DD)$.

Conversely, suppose that $u\in h^p_{\alpha}(\DD)$ for a given $p\in[1,\infty)$. Let
\[
\f_n[|u|^p](x):=\EE^x\left|u(X(\tau_{\DD_n}))\right|^p,\quad n=2,3,...
\]
By the Jensen inequality we have $|u|^p\leq\f_n[|u|^p]$ for every $n$. Therefore, by Fubini theorem, for any $r\in[0,1-1/n)\cup(1+1/n,n)$ we obtain
\[
\|u_r\|^p_p\leq\left\|\f_n[|u|^p]_r\right\|_1
\]
\[
=\int_{\s}\int_{(\overline{\DD}_n)^c}P_{\DD_n}(rx,y)|u(y)|^pdy\sigma(dx)=\int_{(\overline{\DD}_n)^c}\int_{\s}P_{\DD_n}(rx,y)\sigma(dx)|u(y)|^pdy
\]
\[
=\left(\int^{1+1/n}_{1-1/n}+\int^{\infty}_{n}\right)\int_{\s}\int_{\s}P_{\DD_n}(rx,sz)\sigma(dx)|u(sz)|^p\sigma(dz)s^{d-1}ds.
\]
By Lemma~\ref{lem:rotation}, the last term is equal to
\[
\left(\int^{1+1/n}_{1-1/n}+\int^{\infty}_{n}\right)f_n(r,s)\|u_s\|^p_ps^{d-1}ds
\]
\[
\leq\|u\|^p_{h^p}\left(\int^{1+1/n}_{1-1/n}+\int^{\infty}_{n}\right)\int_{\s}P_{\DD_n}(rw,sx)\sigma(dx)s^{d-1}ds=\|u\|^p_{h^p}.
\]
In the last integral, $w\in\s$ is arbitrary. Therefore, $\f_n[|u|^p]$ is finite $\sigma$-a.e. on $r\s$. 
In view of the Harnack inequality (see \cite[Theorem 1]{BogKK}), $\f_n[|u|^p]$ is finite and $\alpha$-harmonic on $\DD_n$ 
for every $n$. As in the proof of Lemma~\ref{lem:majorant1} we conclude that the sequence $\f_n[|u|^p]$ is nondecreasing and 
$\f_{\DD}[|u|^p]=\lim_n\f_n[|u|^p]$. Therefore, by the monotone convergence theorem we obtain 
\[
\|u_r\|^p_p\leq\|\f_{\DD}[|u|^p]_r\|_1\leq\|u\|^p_{h^p}<\infty, 
\]
and from Lemma~\ref{lem:majorant1} it follows that $\f_{\DD}[|u|^p]$ is finite and $\alpha$-harmonic on $\DD$. Taking the supremum over $r$ we obtain 
$\|u\|^p_{h^p}=\|\f_{\DD}[|u|^p]\|_{h^1}$, as desired. 
\end{proof}

\begin{corollary}\label{cor:harmaj}
Let $u$ be $\alpha$-harmonic on $\DD$. Then $u\in h^p_{\alpha}(\DD)$ for a given $p\in[1,\infty)$ if and only if there exists a nonnegative
$\alpha$-harmonic function $v$ on $\DD$ such that $|u|^p\leq v$.
\end{corollary}

\begin{proof}
If such a $v$ exists, then as in the proof of Lemma~\ref{lem:majorant1} we get $\f_{\DD}[|u|^p]\leq v$ and the corollary follows from Lemma~\ref{lem:harmaj}.
\end{proof}

\noindent The next theorem completes the characterization of the spaces $h^p_{\alpha}(\DD)$, $1\leq p\leq\infty$, that we introduced 
in Theorem~\ref{th:hpnorm}.

\begin{theorem}\label{th:spherehp}
Let $u$ be $\alpha$-harmonic on $\DD$.
\begin{enumerate}
	\item[\it 1.] If $u\in h^1_{\alpha}(\DD)$, then $u(x)=\p_{\DD}[\mu](x)+c(1-\Phi(x))$ for some pair $(\mu,c)\in\M(\s)\times\RR$. Furthermore, 
	$\mu$ and $c$ are unique.
	\item[\it 2.] If $u\in h^p_{\alpha}(\DD)$ for a given $p\in(1,\infty]$, then $u(x)=\p_{\DD}[f](x)+c(1-\Phi(x))$ for some pair $(f,c)\in L^p(\s,\sigma)\times\RR$. Furthermore, $f$ and $c$ are unique.
\end{enumerate}
\end{theorem}

\begin{proof}
The first part of the theorem follows immediately from Corollary~\ref{cor:harmaj} and Corollary~\ref{cor:Martinunique1}. 
To prove the second part choose $p\in(1,\infty]$ and suppose that 
$u\in h^p_{\alpha}(\DD)$. Then $u\in h^1_{\alpha}(\DD)$ and by the first part, $u(x)=\mathcal{P}_{\DD}[\mu](x)+c(1-\Phi(x))$ for a unique pair 
$(\mu,c)\in\M(\s)\times\RR$. Hence, it suffices to show that $d\mu=fd\sigma$ and $f\in L^p(\s,\sigma)$. 
Since $\Phi$ is bounded, $c(1-\Phi(\cdot))\in h^p_{\alpha}(\DD)$ and thus $\p_{\DD}[\mu]\in h^p_{\alpha}(\DD)$. 
Therefore, the family $\left\{\p_{\DD}[\mu]_r:r>0\wedge r\neq1\right\}$ is norm-bounded in $L^p(\s,\sigma)$. By Banach-Alaoglu theorem, there is a sequence $r_n$ 
tending to 1 such that $\p_{\DD}[\mu]_{r_n}$ tends weak$^*$ to some $f\in L^p(\s,\sigma)$, i.e., for $q=p/(p-1)$ ($q=1$ when $p=\infty$) and every $g\in L^q(\s,\sigma)$ we have
\[
\int_{\s}g(x)\p_{\DD}[\mu](r_nx)\sigma(dx)\stackrel{n\to\infty}{\longrightarrow}\int_{\s}g(x)f(x)\sigma(dx).
\]
On the other hand, by Lemma~\ref{lem:poissinteg} (v), for every $g\in C(\s)$ we have
\[
\int_{\s}g(x)\p_{\DD}[\mu](r_nx)\sigma(dx)\stackrel{n\to\infty}{\longrightarrow}\int_{\s}g(x)\mu(dx).
\]
Since $C(\s)\subset L^q(\s,\sigma)$, we have $d\mu=fd\sigma$.
\end{proof}

\noindent We will now show that the Fatou theorem holds for $\alpha$-harmonic functions on $\DD$. 
For $y\in\s$ and $\beta>0$ we define the cone $\Gamma_{\beta}(y)$ as
\[
\Gamma_{\beta}(y)=\left\{x\in\DD:|x-y|<(1+\beta)|1-|x||\right\}.
\]
We say that a function $u$ on $\DD$ has a nontangential limit $L$ at $y\in\s$ if, for every $\beta>0$,
\[
\lim_{\Gamma_{\beta}(y)\ni x\to y}u(x)=L.
\]

\begin{theorem}\label{th:Fatou}
Let $\mu\in\M(\s)$ and let $\mu(dx)=f(x)\sigma(dx)+\nu(dx)$ be the Lebesgue decomposition of $\mu$ 
with respect to $\sigma$. Then $\p_{\DD}[\mu]$ has the nontangential limit $f(y)$ at $\sigma$-almost every $y\in\s$.
\end{theorem}

\begin{proof}
For $y\in\s$ and $r>0$ let $K(y,r):=B(y,r)\cap\s$. Define
\[
	\mathcal{L}[\mu](y):=\sup_{r>0}\frac{|\mu|(K(y,r))}{\sigma(K(y,r))}.
\]
Following \cite[Theorem 6.39 and 6.42]{ABR} it is enough to show that for any $\beta>0$ there is a constant $c=c(\alpha,\beta,d)>0$ 
such that for every $y\in\s$ we have
\[
	\sup_{x\in\Gamma_{\beta}(y)}|\p_{\DD}[\mu](x)|\leq c\mathcal{L}[\mu](y).
\]
Fix $y\in\s$ and $x\in\Gamma_{\beta}(y)$. If $|x|\geq2$, then $\p_{\DD}(x,z)\leq c$, where $c=c(\alpha,d)>0$, so 
\[
	|\p_{\DD}[\mu](x)|\leq\int_{\s}\p_{\DD}(x,z)|\mu|(dz)\leq c\|\mu\|\leq c\mathcal{L}[\mu](y).
\]
Suppose that $|x|<2$ and denote $\eta=|x-y|$. We have
\[
	|\p_{\DD}[\mu](x)|\leq\int_{\s}\p_{\DD}(x,z)|\mu|(dz)=\int_{|z-y|<2\eta}\p_{\DD}(x,z)|\mu|(dz)
\]
\[
	+\sum^{\infty}_{k=2}\int_{2^{k-1}\eta<|z-y|<2^k\eta}\p_{\DD}(x,z)|\mu|(dz).
\]
Furthermore, 
\[
	\p_{\DD}(x,z)\leq c_1\frac{||x|-1|^{\alpha-1}}{\left|x-y\right|^{d+\alpha-2}}
	\leq\frac{c_1}{||x|-1|^{d-1}}\leq c_2\eta^{1-d},
\]
where $c_2=c_2(\alpha,\beta,d)>0$. Hence
\[
	\int_{|z-y|<2\eta}\p_{\DD}(x,z)|\mu|(dz)\leq c_2\eta^{1-d}|\mu|(K(y,2\eta))
\]
\[
	\leq c_3\frac{|\mu|(K(y,2\eta))}{\sigma(K(y,2\eta))}\leq c_3\mathcal{L}[\mu](y),
\]
for some $c_3=c_3(\alpha,\beta,d)>0$. For $2^{k-1}\eta<|z-y|<2^k\eta$, $k\geq2$ we have
\[
	|x-z|\geq|z-y|-|y-x|\geq2^{k-1}\eta-\eta\geq2^{k-2}\eta,
\]
and so
\[
	\p_{\DD}(x,z)\leq c_1\frac{||x|-1|^{\alpha-1}}{\left|x-y\right|^{d+\alpha-2}}\leq c_1\frac{\eta^{\alpha-1}}{(2^{k-2}\eta)^{d+\alpha-2}}
=\frac{c_4}{2^{k(d+\alpha-2)}\eta^{d-1}}
\]
\[
	\leq\frac{c_5}{2^{k(\alpha-1)}\sigma(K(y,2^k\eta))},
\]
where $c_5=c_5(\alpha,\beta,d)>0$. Therefore
\[
	\int_{2^{k-1}\eta<|z-y|<2^k\eta}\p_{\DD}(x,z)|\mu|(dz)\leq\frac{c_5|\mu|(K(y,2^k\eta))}{2^{k(\alpha-1)}\sigma(K(y,2^k\eta))}
	\leq\frac{c_5}{2^{k(\alpha-1)}}\mathcal{L}[\mu](y).
\]
Since $\alpha>1$ we have $c_6=\sum^{\infty}_{k=2}2^{-k(\alpha-1)}<\infty$ and hence
\[
|\p_{\DD}[\mu](x)|\leq c_3\mathcal{L}[\mu](y)+\sum^{\infty}_{k=2}\frac{c_5}{2^{k(\alpha-1)}}\mathcal{L}[\mu](y)=(c_3+c_5c_6)\mathcal{L}[\mu](y).
\]
\end{proof}

\begin{corollary}\label{cor:Fatou}
Suppose $u$ is $\alpha$-harmonic and nonnegative on $\DD$. Then $u$ has a nontangential limit at $\sigma$-almost every $y\in\s$.
\end{corollary}

\begin{proof}
The corollary follows from Proposition~\ref{prop:Martin1} and Theorem~\ref{th:Fatou}.
\end{proof}

\noindent At the end of this section we will look at the probabilistic approach to the Hardy spaces discussed in 
\cite{BDL} and \cite{MR1}.
 
\begin{definition}\label{def:hpprob1} For $p\in[1,\infty)$ we define the space $\h^p_{\alpha}(\DD)$ as the family of functions $u$ $\alpha$-harmonic on 
$\DD$ such that
\[
 \|u\|_{\h^p_{\alpha}}:=\sup_{U\subset\subset\DD}\left(\EE^0|u(X_{\tau_U})|^p\right)^{1/p}<\infty.
\]
\end{definition}
\noindent In view of (\ref{det22}) we have $\|u\|_{\h^p_{\alpha}}=(\f_{\DD}[|u|^p](0))^{1/p}$. Then Lemma~\ref{lem:harmaj} implies immediately that $\mathcal{H}^p_{\alpha}(\DD)=h^p_{\alpha}(\DD)$ for all $p\in[1,\infty)$. To identify the norm $\|u\|_{\h^p_{\alpha}}$ we will need the following 
result.

\begin{lemma}\label{lem:majmeasure}
Let $u$ be $\alpha$-harmonic on $\DD$.
\begin{enumerate}
	\item[\it 1.] If $u(x)=\p_{\DD}[\mu](x)+c(1-\Phi(x))$ for some pair $(\mu,c)\in\M(\s)\times\RR$, then 
$\f_{\DD}[u](x)=\p_{\DD}[|\mu|](x)+|c|(1-\Phi(x))$.
	\item[\it 2.] Let $p\in(1,\infty)$. If $u(x)=\p_{\DD}[f](x)+c(1-\Phi(x))$ for some pair $(f,c)\in L^p(\s,\sigma)\times\RR$, then 
$\f_{\DD}[|u|^p](x)=\p_{\DD}[|f|^p](x)+|c|^p(1-\Phi(x))$.
\end{enumerate}
\end{lemma}

\begin{proof}
Let $u(x)=\p_{\DD}[\mu](x)+c(1-\Phi(x))$ for some pair $(\mu,c)\in\M(\s)\times\RR$ and set 
$h(x)=\p_{\DD}[|\mu|](x)+|c|(1-\Phi(x))$. Then $h$ is nonnegative $\alpha$-harmonic on $\DD$ and $|u|\leq h$. 
Hence $\f_{\DD}[u]\leq h$, and by Lemma~\ref{lem:majorant1}, $\f_{\DD}[u]$ is nonnegative $\alpha$-harmonic on $\DD$. 
By Proposition~\ref{prop:Martin1}, $\f_{\DD}[u](x)=\p_{\DD}[\nu](x)+\tilde{c}(1-\Phi(x))$ for some nonnegative 
$\nu\in\M(\s)$ and $\tilde{c}\geq0$. 
As in the proof of \cite[Theorem 1]{MR1} we conclude that $\nu\leq|\mu|$ and $\tilde{c}\leq|c|$. On the other hand, since 
$\f_{\DD}[u]-u\geq0$, by Proposition~\ref{prop:Martin1} and Corollary~\ref{cor:Martinunique1} we have 
$\f_{\DD}[u](x)-u(x)=\p_{\DD}[\nu-\mu](x)+(\tilde{c}-c)(1-\Phi(x))$, so $\nu-\mu\geq0$ and $\tilde{c}-c\geq0$. Applying the 
same argument to $-u$ and taking into account the Hahn decomposition of $\mu$ we obtain $\nu\geq|\mu|$ and $\tilde{c}\geq|c|$. 
This gives the first part.

Let now $p\in(1,\infty)$ and let $u(x)=\p_{\DD}[f](x)+c(1-\Phi(x))$ for some $(f,c)\in L^p(\s,\sigma)\times\RR$. 
Set $h(x)=\p_{\DD}[|f|^p](x)+|c|^p(1-\Phi(x))$. 
Then by the Jensen inequality we have $|u|^p\leq h$ and thus $\f_{\DD}[|u|^p]\leq h$. Using Lemma~\ref{lem:majorant1}, 
Proposition~\ref{prop:Martin1} and the arguments from the previous part we conclude that 
$\f_{\DD}[|u|^p](x)=\p_{\DD}[g](x)+\tilde{c}(1-\Phi(x))$ for some positive function $g\in L^1(\s,\sigma)$ 
with $g\leq|f|^p$ and 
$0\leq\tilde{c}\leq|c|^p$. On the other hand, by the arguments given in the proof of Theorem~\ref{th:hpnorm} we have
\[
\||f|^p\|_1=\lim_{r\to1}\||u_r|^p\|_1\leq\lim_{r\to1}\|\f_{\DD}[|u|^p]_r\|_1=\|g\|_1.
\]
Hence $g=|f|^p$ $\sigma$-a.e. Furthermore, 
\[
|c|^p=\lim_{x\to\infty}|u(x)|^p\leq\lim_{x\to\infty}\f_{\DD}[|u|^p](x)=\tilde{c},
\]
so $\tilde{c}=|c|^p$.
\end{proof}

\noindent By Theorem~\ref{th:hpnorm}, Theorem~\ref{th:spherehp} and Lemma~\ref{lem:majmeasure} we obtain the following complete characterization 
of the spaces $\h^p_{\alpha}(\DD)$.

\begin{theorem}\label{th:hpprob1}
Let $u$ be $\alpha$-harmonic on $\DD$. Then
\begin{enumerate}
	\item[\it 1.] $u\in\h^1_{\alpha}(\DD)$ if and only if $u(x)=\p_{\DD}[\mu](x)+c(1-\Phi(x))$ for some pair $(\mu,c)\in\M(\s)\times\RR$.
	Furthermore, $\mu$ and $c$ are unique and \[\|u\|_{\h^1_{\alpha}}=\Phi(0)\|\mu\|+|c|(1-\Phi(0)).\]
	\item[\it 2.] $u\in\h^p_{\alpha}(\DD)$ for a given $p\in(1,\infty)$ if and only if $u(x)=\p_{\DD}[f](x)+c(1-\Phi(x))$ for some pair 
	$(f,c)\in L^p(\s,\sigma)\times\RR$. Furthermore, $f$ and $c$ are unique and \[\|u\|_{\h^p_{\alpha}}=\left[\Phi(0)\|f\|^p_p+|c|^p(1-\Phi(0))\right]^{1/p}.\]
\end{enumerate}
\end{theorem}

\noindent In particular, for any $p\in[1,\infty)$ and any $\alpha$-harmonic function $u$ on $\DD$ we have
\[
 \left[\Phi(0)\wedge(1-\Phi(0))\right]\|u\|_{h^p}\leq\|u\|_{\h^p_{\alpha}}\leq\|u\|_{h^p}.
\]
We will see in the next section that such a relation does not hold for analogous constructions of Hardy spaces on $\HH$.

\section{$\alpha$-harmonic functions on the complement of the hyperplane}\label{sec4}

We will study the behavior of $\alpha$-harmonic functions on $\HH=\RR^d\setminus\LL$. We will identify 
$\LL=\partial\HH=\left\{x=(x_1,...,x_d)\in\RR^d:x_d=0\right\}$ with the euclidean space $\RR^{d-1}$ and for $x\in\RR^d$ we will denote 
$x=(\x,x_d)$, where $\x\in\RR^{d-1},x_d\in\RR$. 
Let $C(\RR^{d-1})$ denote the space of continuous functions on $\RR^{d-1}$ and let $\M(\RR^{d-1})$ be the space of finite signed Borel 
measures on $\RR^{d-1}$ with the total variation norm $\|\cdot\|$. In this section, we also denote 
by $C_b(\RR^{d-1}),C_0(\RR^{d-1}),C_c(\RR^{d-1})$ the subspaces of $C(\RR^{d-1})$ of functions bounded on $\RR^{d-1}$, vanishing at 
$\infty$ at having a compact support on $\RR^{d-1}$, respectively.
For simplicity, we will adapt the notation of the $L^p$-norm from the previous section, i.e., for a Borel function $f$ on $\RR^{d-1}$ 
and $1\leq p<\infty$ let
\[
\|f\|_p:=\left(\int_{\RR^{d-1}}|f(\x)|^pd\x\right)^{1/p},
\]
and let $\|f\|_{\infty}$ denote the essential supremum norm on $\RR^{d-1}$ with respect to the $(d-1)$-dimensional Lebesgue measure.

We shall return to the question of the hitting probability of $\LL$ for $X_t$.
As we mentioned in the Introduction, the last coordinate of $X_t=(X^1_t,...,X^d_t)$ is a one-dimensional, symmetric $\alpha$-stable 
L\'evy process. Since we consider only $\alpha\in(1,2)$, such a process is pointwise recurrent (see \cite{P2}). Hence
\begin{equation}\label{det32}
\PP^x(T_{\LL}<\infty)=\PP^{x_d}(T_{\left\{0\right\}}<\infty)=1,\quad x\in\HH.
\end{equation}
Here $\PP^{x_d}$ means the distribution of the one-dimensional process $X^d_t$ starting from $x_d$ and $T_{\left\{0\right\}}:=\inf\left\{t>0:X^d_t=0\right\}$. 
We will now calculate the hitting distribution of $\LL$ for $X_t$ (i.e., the $\alpha$-harmonic measure for $\HH$). The formula has been 
proved in \cite{I} only in two-dimensional case. Recent results of \cite{ByMR2} allow us to give the proof for all $d\geq2$.

\begin{proposition}\label{prop:poisker}
The $\alpha$-harmonic measure for $\HH$ has a density with respect to the $(d-1)$-dimensional Lebesgue measure on $\LL$ given by
\[
\p_{\HH}(x,y)=\C_3\frac{|x_d|^{\alpha-1}}{\left|x-y\right|^{d+\alpha-2}},\quad x\in\HH, y\in\LL,
\]
where $\C_3=\pi^{(1-d)/2}\Gamma((\alpha+d)/2-1)/\Gamma((\alpha-1)/2)$.
\end{proposition}

\begin{proof}
We use the methods of \cite[Section 3]{ByMR2}. Let $\Y(t)=(B^d(t),Y(t))$ be a $(d+1)$-dimensional diffusion with independent 
components, where $B^d(t)=(B_1(t),...,B_d(t))$ is the standard Brownian motion in $\RR^d$ and $Y(t)$ is the Bessel process 
with index $-\alpha/2$. Let 
\[
\widetilde{\LL}=\left\{\tilde{x}=(x_1,...,x_{d+1})\in\RR^{d+1}:x_d=0\wedge x_{d+1}=0\right\}
\]
and set $\widetilde{\HH}=\RR^{d+1}\setminus\widetilde{\LL}$. Let $\tau_{\widetilde{\HH}}=\inf\left\{t>0:\Y(t)\notin\widetilde{\HH}\right\}$ 
and define $Z(t)=\sqrt{(B_d(t))^2+(Y(t))^2}$. Then $Z(t)$ is the Bessel process with index 
$\delta=(1-\alpha)/2$. Since $1<\alpha<2$ we have $-1/2<\delta<0$ and hence, for any $a>0$ and $T_0=\inf\left\{t>0:Z(t)=0\right\}$ we have $\PP^a(T_0<\infty)=1$. 
Therefore, for any $\tilde{x}=(x_1,...,x_d,0)\in\HH\times\left\{0\right\}$ we obtain
\[
\PP^{\tilde{x}}(\tau_{\widetilde{\HH}}<\infty)=\PP^{|x_d|}(T_0<\infty)=1.
\]
Moreover, since $0$ is regular for $Z(t)$, every point of $\widetilde{\LL}$ is regular for $\widetilde{\LL}$ with respect to $\Y$. 
Let $x\in\HH$ and set $T^{|x_d|}_0=T_0$ with the starting point $Z_0=|x_d|$. As before we set $\tau_{\HH}=\inf\left\{t>0:X(t)\notin\HH\right\}$. 
Then by \cite[Proposition 3.1, see also Lemma 6.2, Lemma 6.4 and Corollary 6.5]{ByMR2}, for $A\subset\RR^{d-1}$ we have
\[
\PP^x\left(X({\tau_{\HH}})\in A\times\left\{0\right\}\right)=\PP^{\x}\left(B^{d-1}\left(T^{|x_d|}_0\right)\in A\right).
\]
By \cite[p. 75]{BorS} we have
\[
\PP^a(T_0\in dt)=\frac{-2\delta}{a^{2\delta}}\cdot\frac{t^{\delta-1}}{2^{1-\delta}\Gamma(1-\delta)}\exp\left(\frac{-a^2}{2t}\right)dt,\quad a>0.
\]
Since $B^{d-1}$ and $Z$ are independent, we obtain
\[
\PP^{\x}\left(B^{d-1}\left(T^{|x_d|}_0\right)\in A\right)=\int^{\infty}_0\PP^{\x}(B^{d-1}(t)\in A)\PP^{|x_d|}(T_0\in dt)
\]
\[
=\int^{\infty}_0\int_A\frac{1}{(2\pi t)^{\frac{d-1}{2}}}\exp\left(\frac{-|\x-\y|^2}{2t}\right)
\frac{(\alpha-1)t^{\frac{-(\alpha+1)}{2}}}{|x_d|^{1-\alpha}2^{\frac{\alpha+1}{2}}\Gamma\left(\frac{\alpha+1}{2}\right)}\exp\left(\frac{-|x_d|^2}{2t}\right)d\y dt
\]
\[
=\frac{(\alpha-1)|x_d|^{\alpha-1}}{2^{\frac{d+\alpha}{2}}\pi^{\frac{d-1}{2}}\Gamma\left(\frac{\alpha+1}{2}\right)}
\int_A\int^{\infty}_0\left(\frac{1}{t}\right)^{\frac{d+\alpha}{2}}\exp\left(\frac{-\left(|\x-\y|^2+|x_d|^2\right)}{2t}\right)dtd\y.
\]
Because for $a>1$ and $b>0$,
\[
\int^{\infty}_0t^{-a}e^{-b/t}dt=b^{1-a}\int^{\infty}_0s^{a-2}e^{-s}ds=
b^{1-a}\Gamma(a-1),
\]
we obtain
\[
\PP^x\left(X({\tau_{\HH}})\in A\times\left\{0\right\}\right)=
\]
\[
=\frac{(\alpha-1)\Gamma\left(\frac{d+\alpha}{2}-1\right)|x_d|^{\alpha-1}}{2^{\frac{d+\alpha}{2}}\pi^{\frac{d-1}{2}}\Gamma\left(\frac{\alpha+1}{2}\right)}\int_A\left(\frac{2}{|\x-\y|^2+|x_d|^2}\right)^{\frac{d+\alpha-2}{2}}d\y
\]
\[
=\frac{\Gamma\left(\frac{d+\alpha}{2}-1\right)}{\pi^{\frac{d-1}{2}}\Gamma\left(\frac{\alpha-1}{2}\right)}\int_A\frac{|x_d|^{\alpha-1}}
{|x-(\y,0)|^{d+\alpha-2}}d\y.
\]
\end{proof}

\vspace{5mm}

\noindent We will call $\p_{\HH}(x,y)$ the Poisson kernel of $\HH$ for $X_t$. A simple consequence of Proposition~\ref{prop:poisker} is the 
following symmetry property
\begin{equation}\label{det7}
\p_{\HH}((\x,t),y)=\p_{\HH}((\y,t),x),\quad x,y\in\LL,t\in\RR\setminus\left\{0\right\}.
\end{equation}
Let $G_{\HH}(x,y)$ be the Green function of $\HH$ for $X_t$ and let $e_d=(0,...,0,1)$. We define the Martin kernel of $\HH$ for $X_t$ by 
\begin{equation}\label{det10}
M_{\HH}(x,z)=\lim_{\HH\ni y\to z}\frac{G_{\HH}(x,y)}{G_{\HH}(e_d,y)},\quad x\in\HH,z\in\LL\cup\left\{\infty\right\}.
\end{equation}
By \cite[Theorem 2]{BogKK}, the limit in (\ref{det10}) always exists. We will calculate $G_{\HH}$ and $M_{\HH}$ using the methods 
of \cite{BogZ}. The inversion with respect to $\s$ is defined by
\begin{equation}\label{det26}
Tx=\begin{cases}
		x/|x|^2, & \ x\in \RR^d\setminus\left\{0\right\},\\
		\infty, & \ x=0,\\
		0, & \ x=\infty.
		\end{cases}
\end{equation}
This map takes spheres containing 0 onto hyperplanes. Let $\widetilde{T}$ be the inversion with respect to the sphere $S(-e_d,\sqrt{2})$, i.e., 
\begin{equation}\label{det14}
 \widetilde{T}x:=2T(x+e_d)-e_d,
\end{equation}
and let $\HH':=\HH\setminus\left\{-e_d\right\}$. Then we have $\HH'=\widetilde{T}(\DD)$ and $\LL=\widetilde{T}(\s\setminus\left\{-e_d\right\})$. 
By \cite[Theorem 2]{BogZ}, the scaling property and the translation invariance of the symmetric stable processes we obtain
\begin{equation}\label{det15}
G_{\HH'}(x,y)=2^{d-\alpha}|x+e_d|^{\alpha-d}|y+e_d|^{\alpha-d}G_{\DD}(\widetilde{T}x,\widetilde{T}y).
\end{equation}

\begin{proposition}\label{prop:greenmartin}
We have
\begin{equation}\label{det11}
G_{\HH}(x,y)=\frac{\A_{d,\alpha}}{|x-y|^{d-\alpha}}\left[1-\phi\left(\sqrt{1+\frac{4x_dy_d}{|x-y|^2}}\right)\right],\quad x,y\in\HH,
\end{equation}
where $\phi$ is the hitting probability given in (\ref{det3}) and $\A_{d,\alpha}$ is defined in (\ref{det27}). Furthermore,
\begin{equation}\label{det12}
M_{\HH}(x,z)=\frac{\p_{\HH}(x,z)}{\p_{\HH}(e_d,z)}=\frac{|x_d|^{\alpha-1}\left|e_d-z\right|^{d+\alpha-2}}{\left|x-z\right|^{d+\alpha-2}},\quad x\in\HH, z\in\LL,
\end{equation}
and
\begin{equation}\label{det13}
M_{\HH}(x,\infty)=|x_d|^{\alpha-1},\quad x\in\HH.
\end{equation}
\end{proposition}

\begin{proof}
For all $x,y\in\RR^d\setminus\left\{0\right\}$ we have $|Tx-Ty|=|x-y|/(|x||y|)$. Hence we obtain 
\[
 |\widetilde{T}x-\widetilde{T}y|=\frac{2|x-y|}{|x+e_d||y+e_d|},\quad x,y\in\RR^d\setminus\left\{-e_d\right\}.
\]
By (\ref{det15}) we obtain
\[
 G_{\HH'}(x,y)=\A_{d,\alpha}\frac{2^{d-\alpha}|x+e_d|^{\alpha-d}|y+e_d|^{\alpha-d}}
{|\widetilde{T}x-\widetilde{T}y|^{d-\alpha}}\times
\]
\[
 \left[1-\Phi\left(\frac{\widetilde{T}y}{|\widetilde{T}x-\widetilde{T}y|}\left|\widetilde{T}x-\frac{\widetilde{T}y}
{|\widetilde{T}y|^2}\right|\right)\right]=\frac{\A_{d,\alpha}}{|x-y|^{d-\alpha}}\left[1-\Phi\left(N(x,y)\right)\right],
\]
where
\[
 N(x,y)=\frac{\widetilde{T}y}{|\widetilde{T}x-\widetilde{T}y|}\left|\widetilde{T}x-\frac{\widetilde{T}y}
{|\widetilde{T}y|^2}\right|=\frac{\widetilde{T}y}{|\widetilde{T}y||\widetilde{T}x-\widetilde{T}y|}\left||\widetilde{T}y|\widetilde{T}x
-\frac{\widetilde{T}y}{|\widetilde{T}y|}\right|
\]
\[
 =\frac{\widetilde{T}y}{|\widetilde{T}y||\widetilde{T}x-\widetilde{T}y|}
\left(|\widetilde{T}y|^2|\widetilde{T}x|^2-2\langle\widetilde{T}x,\widetilde{T}y\rangle+1\right)^{\frac{1}{2}}
\]
\[
 =\frac{\widetilde{T}y}{|\widetilde{T}y||\widetilde{T}x-\widetilde{T}y|}
\left[|\widetilde{T}x-\widetilde{T}y|^2+\left(1-|\widetilde{T}x|^2\right)\left(1-|\widetilde{T}y|^2\right)\right]^{\frac{1}{2}}
\]
\[
 =\frac{\widetilde{T}y}{|\widetilde{T}y|}
\left[1+\frac{\left(1-|\widetilde{T}x|^2\right)\left(1-|\widetilde{T}y|^2\right)}{|\widetilde{T}x-\widetilde{T}y|^2}\right]^{\frac{1}{2}}.
\]
Furthermore,
\[
 1-|\widetilde{T}x|^2=4\langle T(x+e_d),e_d\rangle-4|T(x+e_d)|^2=\frac{4x_d}{|x+e_d|^2},
\]
so
\[
N(x,y)=\frac{\widetilde{T}y}{|\widetilde{T}y|}\sqrt{1+\frac{4x_dy_d}{|x-y|^2}}.
\]
Since $\left\{-e_d\right\}$ is a polar set, for all $x,y\in\HH'$ we have $G_{\HH}(x,y)=G_{\HH'}(x,y)$, and 
the continuity of $G_{\HH}$ gives (\ref{det11}).

To obtain (\ref{det12}) and (\ref{det13}) we observe first, that $\widetilde{T}e_d=0$. From (\ref{det10}) and (\ref{det15}) 
it follows that
\[
 M_{\HH}(x,z)=2^{d-\alpha}|x+e_d|^{\alpha-d}M_{\DD}(\widetilde{T}x,\widetilde{T}z),\quad x\in\HH',z\in\LL\cup\left\{\infty\right\}.
\]
By (\ref{det8}), for $x\in\HH'$ and $z\in\LL$ we have
\[
 M_{\HH}(x,z)=2^{d-\alpha}|x+e_d|^{\alpha-d}\frac{||\widetilde{T}x|^2-1|^{\alpha-1}}
{|\widetilde{T}x-\widetilde{T}z|^{d+\alpha-2}}=
\frac{|x_d|^{\alpha-1}|e_d-z|^{d+\alpha-2}}{\left|x-z\right|^{d+\alpha-2}}.
\]
Furthermore, since $\widetilde{T}\infty=-e_d$, for $x\in\HH'$ we have
\[
 M_{\HH}(x,\infty)=2^{d-\alpha}|x+e_d|^{\alpha-d}\frac{||\widetilde{T}x|^2-1|^{\alpha-1}}
{|\widetilde{T}x+e_d|^{d+\alpha-2}}=|x_d|^{\alpha-1}.
\]
Finally, the continuity of $M_{\HH}(\cdot,z)$ for all $z\in\LL\cup\left\{\infty\right\}$ gives (\ref{det12}) and (\ref{det13}), so the 
proposition is proved.
\end{proof}

\noindent We remark that, in opposite to the previous case of $\DD$, the Poisson kernel and the Martin kernel of $\HH$ are no longer the same objects. 
In the next part of this section, for $x\in\HH$ and $\y\in\RR^{d-1}$ we will use the notation $\p_{\HH}(x,\y):=\p_{\HH}(x,y)$ where $y=(\y,0)\in\LL$. 
Analogously we define $M_{\HH}(x,\y)$.

\begin{proposition}\label{prop:Martin2}
For every nonnegative measure $\mu\in\M(\RR^{d-1})$ and every constant $c\geq0$ the function $u$ given by
\begin{equation}\label{det18}
 u(x)=\int_{\RR^{d-1}}M_{\HH}(x,\y)\mu(d\y)+c|x_d|^{\alpha-1},\quad x\in\HH,
\end{equation}
is $\alpha$-harmonic on $\HH$. Conversely, if $u$ is nonnegative and $\alpha$-harmonic on $\HH$ then there exists a unique nonnegative measure 
$\mu\in\M(\RR^{d-1})$ and a unique constant $c\geq0$ satisfying (\ref{det18}).
\end{proposition}

\begin{proof}
We use similar arguments as in the proof of Proposition~\ref{prop:Martin1}. However, it could be slightly more complicated to see directly from (\ref{det11}) 
that $\int_{\HH}G_{\HH}(e_d,y)dy=\infty$ (the accessibility of the point at infinity). We may avoid this difficulty by the fact that $\infty$ is 
accessible from $\HH$ if and only if 0 is accessible from $T(\HH)=\HH$ (see \cite[(82)]{BogKK}).
\end{proof}

For $\mu\in\M(\RR^{d-1})$ and $f\in L^p(\RR^{d-1})$, $1\leq p\leq\infty$, we define the Poisson integrals of $\mu$ and $f$ on $\HH$ as
\[
 \p_{\HH}[\mu](x)=\int_{\RR^{d-1}}\p_{\HH}(x,\y)\mu(d\y),\quad \p_{\HH}[f](x)=\int_{\RR^{d-1}}\p_{\HH}(x,\y)f(\y)d\y.
\]
Analogously we define the Martin integral $M_{\HH}[\mu]$ for $\mu\in\M(\RR^{d-1})$. 
Proposition~\ref{prop:Martin2}, (\ref{det12}) and the Hahn decomposition for signed measures imply that $M_{\HH}[\mu]$, $\p_{\HH}[\mu]$ 
and $\p_{\HH}[f]$ are $\alpha$-harmonic on $\HH$ for every $\mu\in\M(\RR^{d-1})$ and $f\in L^p(\RR^{d-1})$, $1\leq p\leq\infty$.
Furthermore, by (\ref{det12}) we get that $\p_{\HH}[\mu]$ is well-defined and $\alpha$-harmonic on $\HH$ for a measure $\mu$ with not necessarily 
finite variation, but verifying
\begin{equation}\label{det37}
\int_{\RR^{d-1}}\frac{|\mu|(d\x)}{|e_d-(\x,0)|^{d+\alpha-2}}<\infty.
\end{equation}
To simplify the notation, we let $\omega_{\alpha}:=\omega^{e_d}(\cdot,\HH)$ be the $\alpha$-harmonic measure for $\HH$ with the starting point $e_d$. 
Clearly, by Proposition~\ref{prop:poisker}, $\omega_{\alpha}$ is a probability measure on $\RR^{d-1}$ given by
\begin{equation}\label{det21}
 \omega_{\alpha}(d\x)=\mathcal{P}_{\HH}(e_d,\x)d\x=\frac{\C_3d\x}{|e_d-(\x,0)|^{d+\alpha-2}}.
\end{equation}
In view of (\ref{det37}), for any $f\in L^p(\RR^{d-1},\omega_{\alpha})$, $1\leq p<\infty$, 
$\p_{\HH}[f]$ is well-defined and $\alpha$-harmonic on $\HH$. We remark that $L^p(\RR^{d-1},\omega_{\alpha})$ is essentially 
bigger than $L^p(\RR^{d-1})$ for every $p\in[1,\infty)$. Let $\|\cdot\|_{p,\alpha}$ denote the norm associated with 
the space $L^p(\RR^{d-1},\omega_{\alpha})$. We also adapt to the present case the notation introduced in (\ref{det33}), i.e., 
for $t\in\RR\setminus\left\{0\right\}$ and a function $u$ on $\HH$ we define the function $u_t$ on $\RR^{d-1}$ by
\begin{equation}\label{det35}
u_t(\x):=u(\x,t),\quad\x\in\RR^{d-1}.
\end{equation}
The next 3 lemmas characterize the behavior of the Poisson and Martin integrals on $\HH$.

\begin{lemma}\label{lem:poissinteg2}
Poisson integrals on $\HH$ have the following properties:
\begin{enumerate}[\upshape (i)]
 \item If $\mu\in\M(\RR^{d-1})$, then $\|\p_{\HH}[\mu]_t\|_1\leq\|\mu\|$ for every $t$.
 \item If $1\leq p\leq\infty$ and $f\in L^p(\RR^{d-1})$, then $\|\p_{\HH}[f]_t\|_p\leq\|f\|_p$ for every $t$.
 \item If $f\in C_0(\RR^{d-1})$, then $\|\p_{\HH}[f]_t-f\|_{\infty}\to0$ as $t\to0$.
 \item If $1\leq p<\infty$ and $f\in L^p(\RR^{d-1})$, then $\|\p_{\HH}[f]_t-f\|_p\to0$ as $t\to0$.
 \item If $\mu\in\M(\RR^{d-1})$, then $\p_{\HH}[\mu]_t\to\mu$ weak$^*$ in $\M(\RR^{d-1})$ as $t\to0$.
 \item If $f\in L^{\infty}(\RR^{d-1})$, then $\p_{\HH}[f]_t\to f$ weak$^*$ in $L^{\infty}(\RR^{d-1})$ as $t\to0$.
\end{enumerate}
\end{lemma}

\begin{proof}
Using the property (\ref{det7}) and the Jensen inequality, we follow the proofs of the classical counterparts in \cite[Theorems 7.4, 7.6, 7.8 and 7.10]{ABR}.
\end{proof}

\begin{lemma}\label{lem:Martinweak}
Let $u(x)=M_{\HH}[\mu](x)+c|x_d|^{\alpha-1}$ for some $(\mu,c)\in\M(\RR^{d-1})\times\RR$. Then the family of measures
\[
\mu^{\alpha}_t(d\x):=u_t(\x)\omega_{\alpha}(d\x),\quad 0<|t|<\varepsilon,
\]
is norm-bounded in $\M(\RR^{d-1})$ for every $\varepsilon>0$. Furthermore, 
$\mu^{\alpha}_t\to\mu$ weakly as $t\to0$.
\end{lemma}

\begin{proof}
Let $\varepsilon>0$. First we will show that there exists a constant $c_1>0$ depending only on $\varepsilon,\alpha$ and $d$ such that
for every $\y\in\RR^{d-1}$ and $0<|t|<\varepsilon$ we have 
\begin{equation}\label{det20}
 \int_{\RR^{d-1}}M_{\HH}((\x,t),\y)\omega_{\alpha}(d\x)\leq c_1.
\end{equation}
For any $t\neq0$ the left-hand side of (\ref{det20}) is equal to
\[
 \int_{\RR^{d-1}}\frac{|t|^{\alpha-1}|(\y,0)-e_d|^{d+\alpha-2}}{|(\x,t)-(\y,0)|^{d+\alpha-2}}\omega_{\alpha}(d\x)
\]
\[
\leq |t|^{\alpha-1}2^{d+\alpha-2}\int_{\RR^{d-1}}\frac{(|(\y,0)-(\x,t)|\vee|(\x,t)-e_d|)^{d+\alpha-2}}{|(\x,t)-(\y,0)|^{d+\alpha-2}}\omega_{\alpha}(d\x)
\]
\[
\leq |t|^{\alpha-1}2^{d+\alpha-2}\left(1+\int_{\RR^{d-1}}\frac{|(\x,t)-e_d|^{d+\alpha-2}}{|(\x,t)-(\y,0)|^{d+\alpha-2}}\omega_{\alpha}(d\x)\right).
\]
There is a constant $c_2>1$ depending only on $\varepsilon$ such that for every $\x\in\RR^{d-1}$ and $|t|<\varepsilon$ we have 
$|(\x,t)-e_d|\leq c_2|(\x,0)-e_d|$. Hence, for $0<|t|<\varepsilon$ we obtain
\[
\int_{\RR^{d-1}}\frac{|(\x,t)-e_d|^{d+\alpha-2}}{|(\x,t)-(\y,0)|^{d+\alpha-2}}\omega_{\alpha}(d\x)\leq
\int_{\RR^{d-1}}\frac{\C_3c^{d+\alpha-2}_2}{|(\x,t)-(\y,0)|^{d+\alpha-2}}d\x
\]
\[
=\int_{\RR^{d-1}}\frac{\C_3c^{d+\alpha-2}_2}{|(\x,0)-(\y,t)|^{d+\alpha-2}}d\x=|t|^{1-\alpha}c^{d+\alpha-2}_2.
\]
Therefore (\ref{det20}) follows with $c_1=2^{d+\alpha-2}(\varepsilon^{\alpha-1}+c^{d+\alpha-2}_2)$. For $0<|t|<\varepsilon$ we now have
\[
\|\mu^{\alpha}_t\|=\int_{\RR^{d-1}}|M_{\HH}[\mu]_t(\x)+c|t|^{\alpha-1}|\omega_{\alpha}(d\x)\leq\int_{\RR^{d-1}}M_{\HH}[|\mu|]_t(\x)\omega_{\alpha}(d\x)
+|c||t|^{\alpha-1}
\]
\[
=\int_{\RR^{d-1}}\int_{\RR^{d-1}}M_{\HH}((\x,t),\y)\omega_{\alpha}(d\x)|\mu|(d\y)+|c||t|^{\alpha-1}\leq c_1\|\mu\|+|c|\varepsilon^{\alpha-1}.
\]
To prove the second part of the lemma, choose $g\in C_b(\RR^{d-1})$. We have
\[
 \int_{\RR^{d-1}}g(\x)\mu^{\alpha}_t(d\x)=\int_{\RR^{d-1}}\int_{\RR^{d-1}}g(\x)M_{\HH}((\x,t),\y)
\omega_{\alpha}(d\x)\mu(d\y)+c|t|^{\alpha-1}\p_{\HH}[g](e_d).
\]
Obviously, $c|t|^{\alpha-1}\p_{\HH}[g](e_d)$ vanishes as $t\to0$. Furthermore, by (\ref{det20}),
\[
 \left|\int_{\RR^{d-1}}g(\x)M_{\HH}((\x,t),\y)\omega_{\alpha}(d\x)\right|\leq c_1\|g\|_{\infty}.
\]
By (\ref{det12}) and (\ref{det7})
\[
 \int_{\RR^{d-1}}g(\x)M_{\HH}((\x,t),\y)\omega_{\alpha}(d\x)=
\int_{\RR^{d-1}}g(\x)\frac{\p_{\HH}((\x,t),\y)}{\p_{\HH}(e_d,\y)}\p_{\HH}(e_d,\x)d\x
\]
\[
 =\int_{\RR^{d-1}}\frac{\p_{\HH}((\y,t),\x)}{\p_{\HH}(e_d,\y)}g(\x)\p_{\HH}(e_d,\x)d\x
=\frac{\p_{\HH}[g\p_{\HH}(e_d,\cdot)]_t(\y)}{\p_{\HH}(e_d,\y)},
\]
and since $g\p_{\HH}(e_d,\cdot)\in C_0(\RR^{d-1})$, by Lemma~\ref{lem:poissinteg2} (iii) we have 
\[
 \frac{\p_{\HH}[g\p_{\HH}(e_d,\cdot)]_t(\y)}{\p_{\HH}(e_d,\y)}\stackrel{t\to0}{\longrightarrow}g(\y),\quad\y\in\RR^{d-1}.
\]
Hence, by the dominated convergence theorem we obtain
\[
\int_{\RR^{d-1}}\int_{\RR^{d-1}}g(\x)M_{\HH}((\x,t),\y)\omega_{\alpha}(d\x)\mu(d\y)\stackrel{t\to0}{\longrightarrow}
\int_{\RR^{d-1}}g(\y)\mu(d\y),
\]
so the lemma is proved.
\end{proof}

\begin{corollary}\label{cor:Martinunique2}
Let $u$ be $\alpha$-harmonic on $\HH$. Then $u(x)=M_{\HH}[\mu](x)+c|x_d|^{\alpha-1}$ for some pair $(\mu,c)\in\M(\RR^{d-1})\times\RR$ 
if and only if there exists a nonnegative $\alpha$-harmonic function $v$ on $\HH$ such that $|u|\leq v$. Furthermore, $\mu$ and $c$ are unique.
\end{corollary}

\begin{proof}
In view of Proposition~\ref{prop:Martin2} and Lemma~\ref{lem:Martinweak} the proof is similar as in the case of Corollary~\ref{cor:Martinunique1}.
\end{proof}

\begin{lemma}\label{lem:strongconv}
Let $f\in L^p(\RR^{d-1},\omega_{\alpha})$ for a given $p\in[1,\infty)$. Then $\|\p_{\HH}[f]_t-f\|_{p,\alpha}\to0$ as $t\to0$.
\end{lemma}

\begin{proof}
Fix $f\in L^p(\RR^{d-1},\omega_{\alpha})$ and $\varepsilon>0$. Choose $g\in C_c(\RR^{d-1})$ such that 
$\|f-g\|_{p,\alpha}<\varepsilon$. Then we have
\[
 \|\p_{\HH}[f]_t-f\|_{p,\alpha}\leq\|\p_{\HH}[f]_t-\p_{\HH}[g]_t\|_{p,\alpha}+\|\p_{\HH}[g]_t-g\|_{p,\alpha}+\varepsilon.
\]
By Lemma~\ref{lem:poissinteg2} (iii), $\p_{\HH}[g]_t\to g$ uniformly as $t\to0$, so $\|\p_{\HH}[g]_t-g\|_{p,\alpha}<\varepsilon$ 
for $|t|$ sufficiently small. Furthermore, by the Jensen inequality,
\[
 \|\p_{\HH}[f]_t-\p_{\HH}[g]_t\|^p_{p,\alpha}=\int_{\RR^{d-1}}|\p_{\HH}[f]_t(\x)-\p_{\HH}[g]_t(\x)|^p\omega_{\alpha}(d\x)
\]
\[
 \leq\int_{\RR^{d-1}}\int_{\RR^{d-1}}\p_{\HH}((\x,t),\y)|f(\y)-g(\y)|^pd\y\omega_{\alpha}(d\x).
\]
By (\ref{det7}), (\ref{det12}) and Fubini theorem, the last term above is equal to
\[
\int_{\RR^{d-1}}\int_{\RR^{d-1}}M_{\HH}((\x,t),\y)\omega_{\alpha}(d\x)|f(\y)-g(\y)|^p\omega_{\alpha}(d\y),
\]
and by (\ref{det20}), for $0<|t|<1$ we get 
\[
\|\p_{\HH}[f]_t-\p_{\HH}[g]_t\|^p_{p,\alpha}\leq c\|g-f\|^p_{p,\alpha}\leq c\varepsilon^p,
\]
where $c$ depends only on $d$ and $\alpha$. Since $\varepsilon$ was arbitrary, we conclude that $\|\p_{\HH}[f]_t-f\|_{p,\alpha}\to0$ when 
$t\to0$, as desired.
\end{proof}

\noindent We will now describe the corresponding Hardy spaces on $\HH$. In opposite to the previous section, we start with 
the probabilistic definition of $\h^p_{\alpha}(\HH)$.

\begin{definition}\label{def:hpprob2}
For $p\in[1,\infty)$ we define the space $\h^p_{\alpha}(\HH)$ as the family of functions $u$ $\alpha$-harmonic on 
$\HH$ such that
\[
 \|u\|_{\h^p_{\alpha}}:=\sup_{U\subset\subset\HH}\left(\EE^{e_d}|u(X_{\tau_U})|^p\right)^{1/p}<\infty.
\]
\end{definition}

\noindent In view of (\ref{det22}), $\|u\|_{\h^p_{\alpha}}=(\f_{\HH}[|u|^p](e_d))^{1/p}$. By the Jensen inequality, 
for $1\leq p<q\leq\infty$ we have $\|u\|_{\h^p_{\alpha}}\leq\|u\|_{\h^q_{\alpha}}$ and hence 
$\h^q_{\alpha}(\HH)\subset\h^p_{\alpha}(\HH)$. The following is a counterpart of Lemma~\ref{lem:majmeasure}.

\begin{lemma}\label{lem:majorantform} 
Let $u$ be $\alpha$-harmonic on $\HH$.
\begin{enumerate}
\item[\it 1.] If $u(x)=M_{\HH}[\mu](x)+c|x_d|^{\alpha-1}$ for some $(\mu,c)\in\M(\RR^{d-1})\times\RR$, then
 $\f_{\HH}[u](x)=M_{\HH}[|\mu|](x)+|c||x_d|^{\alpha-1}$.
\item[\it 2.] Let $p\in[1,\infty)$. If $u=\p_{\HH}[f]$ for some $f\in L^p(\RR^{d-1},\omega_{\alpha})$, then 
$\f_{\HH}[|u|^p]=\p_{\HH}[|f|^p]$.
\end{enumerate}
\end{lemma}

\begin{proof}
In view of Proposition~\ref{prop:Martin2}, Corollary~\ref{cor:Martinunique2} and Lemma~\ref{lem:strongconv}, the proof is similar as in the case 
of Lemma~\ref{lem:majmeasure}.
\end{proof}

\noindent The next theorem fully characterizes the spaces $\h^p_{\alpha}(\HH)$ in terms of the Martin and Poisson integrals.

\begin{theorem}\label{th:widerhp}
Let $u$ be $\alpha$-harmonic on $\HH$.
\begin{enumerate}
\item[\it 1.] $u\in\h^1_{\alpha}(\HH)$ if and only if $u(x)=M_{\HH}[\mu](x)+c|x_d|^{\alpha-1}$ for some pair $(\mu,c)\in\M(\RR^{d-1})\times\RR$.
Furthermore, $\mu$ and $c$ are unique and $\|u\|_{\h^1_{\alpha}}=\|\mu\|+|c|$.
\item[\it 2.] $u\in\h^p_{\alpha}(\HH)$ for a given $p\in(1,\infty)$ if and only if $u=\p_{\HH}[f]$ for some function $f\in L^p(\RR^{d-1},\omega_{\alpha})$. 
Furthermore, $f$ is unique and $\|u\|_{\h^p_{\alpha}}=\|f\|_{p,\alpha}$.
\end{enumerate}
\end{theorem}

\begin{proof}
If $u(x)=M_{\HH}[\mu](x)+c|x_d|^{\alpha-1}$ for some $(\mu,c)\in\M(\RR^{d-1})\times\RR$, then by Lemma~\ref{lem:majorantform} we have 
$\f_{\HH}[u](x)=M_{\HH}[|\mu|](x)+|c||x_d|^{\alpha-1}$. Hence $\|u\|_{\h^1_{\alpha}}=\f_{\HH}[u](e_d)=\|\mu\|+|c|$ and 
$u\in\h^1_{\alpha}(\HH)$. Conversely, if $u\in\h^1_{\alpha}(\HH)$ then by Lemma~\ref{lem:majorant1}, 
$\f_{\HH}[u]$ is finite and $\alpha$-harmonic on $\HH$ and $|u|\leq\f_{\HH}[u]$. By Corollary~\ref{cor:Martinunique2} we have 
$u(x)=M_{\HH}[\mu](x)+c|x_d|^{\alpha-1}$ for a unique pair $(\mu,c)\in\M(\RR^{d-1})\times\RR$. This proves the first part.

Let now $p\in(1,\infty)$. If $u=\p_{\HH}[f]$ for some $f\in L^p(\RR^{d-1},\omega_{\alpha})$, then by Lemma~\ref{lem:majorantform} we have
$\f_{\HH}[|u|^p]=\p_{\HH}[|f|^p]$. Hence $\|u\|_{\h^p_{\alpha}}=(\f_{\HH}[|u|^p](e_d))^{1/p}=\|f\|_{p,\alpha}$ and
$u\in\h^p_{\alpha}(\HH)$. Conversely, suppose that $u\in\h^p_{\alpha}(\HH)$. Then by Lemma~\ref{lem:majorant1}, $\f_{\HH}[|u|^p]$ is finite 
and $\alpha$-harmonic on $\HH$ and $|u|^p\leq\f_{\HH}[|u|^p]$. 
Since $\h^p_{\alpha}(\HH)\subset\h^1_{\alpha}(\HH)$, by the first part of the theorem we have $u(x)=M_{\HH}[\mu](x)+c_1|x_d|^{\alpha-1}$ 
for a unique pair $(\mu,c_1)\in\M(\RR^{d-1})\times\RR$. 
Because $|x_d|^{\alpha-1}$ is bounded in the neighborhood of $\LL$, there exists 
a nonnegative $\alpha$-harmonic function $v$ on $\HH$ such that $|M_{\HH}[\mu](x)|^p\leq v(x)$ for every $x\in\RR^d$ such that $0<|x_d|<1/2$. 
Furthermore, by Proposition~\ref{prop:Martin2} we have $v(x)=M_{\HH}[\nu](x)+c_2|x_d|^{\alpha-1}$ 
for a unique nonnegative measure $\nu\in\M(\RR^{d-1})$ and a constant $c_2\geq0$. Hence, for $0<|t|<1/2$ 
we obtain 
\[
 \int_{\RR^{d-1}}|M_{\HH}[\mu]_t(\x)|^p\omega_{\alpha}(d\x)\leq\int_{\RR^{d-1}}M_{\HH}[\nu]_t(\x)\omega_{\alpha}(d\x)+c_2(1/2)^{\alpha-1}.
\]
By Lemma~\ref{lem:Martinweak} we conclude, that the family $M_{\HH}[\mu]_t$, $0<|t|<1/2$ is norm-bounded in 
$L^p(\RR^{d-1},\omega_{\alpha})$. By Banach-Alaoglu 
theorem there is a sequence $t_n$ tending to 0 such that $M_{\HH}[\mu]_{t_n}$ tends weak$^*$ to some function 
$f\in L^p(\RR^{d-1},\omega_{\alpha})$, i.e., for $q=p/(p-1)$ and every $g\in L^q(\RR^{d-1},\omega_{\alpha})$ we have
\[
\int_{\RR^{d-1}}g(\x)M_{\HH}[\mu]_{t_n}(\x)\omega_{\alpha}(d\x)\stackrel{n\to\infty}{\longrightarrow}\int_{\RR^{d-1}}g(\x)f(\x)\omega_{\alpha}(d\x).
\]
On the other hand, by Lemma~\ref{lem:Martinweak}, for any $g\in C_b(\RR^{d-1})$ 
\[
\int_{\RR^{d-1}}g(\x)M_{\HH}[\mu]_{t_n}(\x)\omega_{\alpha}(d\x)\stackrel{n\to\infty}{\longrightarrow}\int_{\RR^{d-1}}g(\x)\mu(d\x).
\]
Since $C_b(\RR^{d-1})\subset L^q(\RR^{d-1},\omega_{\alpha})$ we have $\mu(d\x)=f(\x)\omega_{\alpha}(d\x)$ 
and $M_{\HH}[\mu]=\p_{\HH}[f]$. Therefore $u(x)=\p_{\HH}[f](x)+c_1|x_d|^{\alpha-1}$. Because both $|u|^p$ 
and $|\p_{\HH}[f]|^p$ have nonnegative 
$\alpha$-harmonic majorants, we conclude that also $|c_1|^p|x_d|^{p(\alpha-1)}$ has an $\alpha$-harmonic majorant. By Proposition~\ref{prop:Martin2} we have 
$|c_1|^p|x_d|^{p(\alpha-1)}\leq M_{\HH}[\tilde{\nu}](x)+c_3|x_d|^{\alpha-1}$ for a positive measure $\tilde{\nu}\in\M(\RR^{d-1})$ 
and a constant $c_3\geq0$. Taking $x_n=ne_d$, $n=1,2,...$ we obtain
\[
 |c_1|^pn^{p(\alpha-1)}\leq n^{\alpha-1}\int_{\RR^{d-1}}\left(\frac{|e_d-(\y,0)|}{|ne_d-(\y,0)|}\right)^{d+\alpha-2}\tilde{\nu}(d\y)
+c_3n^{\alpha-1}
\]
\[
 \leq n^{\alpha-1}\|\tilde{\nu}\|+c_3n^{\alpha-1}=c_4n^{\alpha-1}.
\]
Since the above estimate is false for $c_1\neq0$ and sufficiently big $n$, we have $c_1=0$ and $u=\p_{\HH}[f]$, as desired.
\end{proof}

\noindent We will now focus on the analytic case. We recall that in this section, the notation $u_t$ is given by (\ref{det35}).

\begin{definition}\label{def:hpan2}
For $p\in [1,\infty]$ we define the space $h^p_{\alpha}(\HH)$ as the family of functions $u$ $\alpha$-harmonic 
on $\HH$ such that
\[
\left\|u\right\|_{h^p}:=\sup_{t\in\RR\setminus\left\{0\right\}}\|u_t\|_p<\infty.
\]
\end{definition}

\noindent We remark that the spaces $h^p_{\alpha}(\HH)$ are slightly more difficult to study than the spaces $h^p_{\alpha}(\DD)$ since 
$\partial\HH=\LL$ is not compact. For $\varepsilon>0$ we define 
\[
\HH_{\varepsilon}:=\left\{x=(x_1,...,x_d)\in\RR^d: |x_d|>\varepsilon\right\}.
\]
We have $\overline{\HH}_{\varepsilon}\subset\HH$ for every $\varepsilon>0$, 
$\HH_{\varepsilon_1}\subset\HH_{\varepsilon_2}$ for $\varepsilon_1>\varepsilon_2>0$ and 
$\bigcup_{\varepsilon>0}\HH_{\varepsilon}=\HH$. $\HH_{\varepsilon}$ will play the role of the sets $\DD_n$ from the previous sections. 
However, since $\HH_{\varepsilon}$ is unbounded, not every $\alpha$-harmonic function on $\HH$ satisfies the mean 
value property (\ref{det30}) with $U=\HH_{\varepsilon}$. A simple counterexample is $|x_d|^{\alpha-1}$. We solve partially this problem 
in the next lemma by giving a sufficient (but not necessary) condition for the property to be verified.

\begin{lemma}\label{lem:meanvalue}
Let $u$ be $\alpha$-harmonic on $\HH$ and let $\varepsilon>0$ be fixed. Suppose that $u$ is bounded on $\HH_{\varepsilon}$ and 
$\EE^{x_0}\left|u\left(X(\tau_{\HH_{\varepsilon}})\right)\right|<\infty$ for some $x_0\in\HH_{\varepsilon}$. Then 
\begin{equation}\label{det19}
 u(x)=\EE^{x}u\left(X(\tau_{\HH_{\varepsilon}})\right),\quad x\in\HH_{\varepsilon}.
\end{equation}
\end{lemma}

\begin{proof}
For $n=1,2,...$ set $U_n:=B(0,n)\cap\HH_{\varepsilon}$. Then for every $n>\varepsilon$ we have 
$\tau_{U_n}\leq\tau_{\HH_{\varepsilon}}$ a.s. and 
\[
 u(x)=\EE^x\left[u\left(X(\tau_{U_n})\right);\tau_{U_n}=\tau_{\HH_{\varepsilon}}\right]+
\EE^x\left[u\left(X(\tau_{U_n})\right);\tau_{U_n}<\tau_{\HH_{\varepsilon}}\right].
\]
Since $u$ is bounded on $\HH_{\varepsilon}$, there is a constant $c>0$ independent of $x$ and $n$ such that
\[
 \EE^x\left[\left|u\left(X(\tau_{U_n})\right)\right|;\tau_{U_n}<\tau_{\HH_{\varepsilon}}\right]\leq c\PP^x(\tau_{U_n}<\tau_{\HH_{\varepsilon}}).
\]
For fixed $x\in\HH_{\varepsilon}$ and $n>|x|$ we have $\PP^x(\tau_{U_n}<\tau_{\HH_{\varepsilon}})\leq\PP^x(\tau_{B(0,n)}\leq\tau_{\HH_{\varepsilon}})$, so
\[
 \lim_{n\to\infty}\PP^x(\tau_{U_n}<\tau_{\HH_{\varepsilon}})\leq\lim_{n\to\infty}\PP^x(\tau_{B(0,n)}\leq\tau_{\HH_{\varepsilon}})=
\PP^x(\bigcap_n\left\{\tau_{B(0,n)}\leq\tau_{\HH_{\varepsilon}}\right\}).
\]
Since 
\[
 \left\{\tau_{\HH_{\varepsilon}}<\infty\right\}\cap\bigcap_n\left\{\tau_{B(0,n)}\leq\tau_{\HH_{\varepsilon}}\right\}=\emptyset
\]
and $\PP^x(\tau_{\HH_{\varepsilon}}<\infty)=1$, we have 
\[
 \lim_{n\to\infty}|\EE^x[u\left(X(\tau_{U_n})\right);\tau_{U_n}<\tau_{\HH_{\varepsilon}}]|\leq 
c\lim_{n\to\infty}\PP^x(\tau_{U_n}<\tau_{\HH_{\varepsilon}})=0.
\]
On the other hand, 
\[
 \EE^x[u\left(X(\tau_{U_n})\right);\tau_{U_n}=\tau_{\HH_{\varepsilon}}]=
\int_{\left(\overline{\HH}_{\varepsilon}\right)^c}P_{U_n}(x,y)u(y)dy.
\]
Since $G_{U_n}(x,y)\nearrow G_{\HH}(x,y)$ as $n\to\infty$, from (\ref{det24}) and the monotone convergence theorem we have $P_{U_n}(x,y)\nearrow P_{\HH_{\varepsilon}}(x,y)$ for every $y\in(\overline{\HH}_{\varepsilon})^c$. Because $\EE^{x_0}\left|u\left(X(\tau_{\HH_{\varepsilon}})\right)\right|<\infty$, by the Harnack inequality (see \cite[Theorem 1]{BogKK}), $\EE^x\left|u\left(X(\tau_{\HH_{\varepsilon}})\right)\right|<\infty$ for every $x$, and the lemma follows from the dominated convergence theorem.
\end{proof}

\begin{lemma}\label{lem:Hpbound}
Suppose that $u\in h^p_{\alpha}(\HH)$ for some $p\in[1,\infty]$. Then $u$ is bounded on $\HH_{\varepsilon}$ for every $\varepsilon>0$.
\end{lemma}

\begin{proof}
Since for $p=\infty$ the result is obvious, we assume that $1\leq p<\infty$. Let $u\in h^p_{\alpha}(\HH)$ and let 
$\varepsilon>0$, $x_0\in\HH_{\varepsilon}$ be fixed. Set $r=\varepsilon/3$ and 
\[
\tau_0=\inf\left\{t>0:X(t)\notin B(x_0,r)\right\},
\quad\tau_x=\inf\left\{t>0:X(t)\notin B(x,2r)\right\}.
\]
For $x\in B(x_0,r)$ let $f_1(x):=\EE^x\left|u(X(\tau_0))\right|^p$, $f_2(x)=\EE^x\left|u(X(\tau_x))\right|^p$. Since 
$B(x_0,r)\subset B(x,2r)$ for all $x\in B(x_0,r)$, we have $\tau_0\leq\tau_x$. Moreover, $B(x,2r)\subset\HH$. By 
the strong Markov property and the Jensen's inequality we have
\[
f_1(x)=\EE^x\left|\EE^{X(\tau_0)}u(X(\tau_x))\right|^p
\leq\EE^x\left[\EE^{X(\tau_0)}\left|u(X(\tau_x))\right|^p\right]=f_2(x).
\]
Furthermore,
\[
\int_{B(x_0,r)}f_2(x)dx=\int_{B(x_0,r)}\int_{2r<|y-x|<3r}P_{B(x,2r)}(x,y)|u(y)|^pdydx
\]
\[
+\int_{B(x_0,r)}\int_{|y-x|>3r}P_{B(x,2r)}(x,y)|u(y)|^pdydx=I_1+I_2.
\]
By (\ref{det23}) we have
\[
I_1=\C_1\int_{B(x_0,r)}\int_{2r<|y-x|<3r}\left(\frac{4r^2}{|x-y|^2-4r^2}\right)^{\alpha/2}\frac{|u(y)|^p}{|x-y|^d}dydx
\]
\[
=\int_{B(x_0,r)}\int_{B(x_0,4r)}P_{B(y,2r)}(y,x)|u(y)|^p1_{\left\{2r<|y-x|<3r\right\}}(x,y)dydx
\]
\[
=\int_{B(x_0,4r)}\int_{B(x_0,r)}P_{B(y,2r)}(y,x)1_{\left\{2r<|y-x|<3r\right\}}(x,y)dx|u(y)|^pdy
\]
\[
\leq\int_{B(x_0,4r)}|u(y)|^pdy\leq\int^{x^d_0+4r}_{x^d_0-4r}\int_{\RR^{d-1}}|u(\y,t)|^pd\y dt\leq8r\left\|u\right\|^p_{h^p},
\]
where $x^d_0$ is the last coordinate of $x_0$. On the other hand,
\[
I_2=\C_1\int_{B(x_0,r)}\int_{|y-x|>3r}\left(\frac{4r^2}{|x-y|^2-4r^2}\right)^{\alpha/2}\frac{|u(y)|^p}{|x-y|^d}dydx
\]
\[
\leq \C_1\int_{B(x_0,r)}\int_{|y-x|>3r}\frac{|u(y)|^p}{|x-y|^d}dydx.
\]
For fixed $x\in B(x_0,r)$ we set
\[
 A_0=\left\{y:x_d-4r<y_d<x_d+4r\right\}\cap\left\{y:|x-y|>3r\right\},
\]
\[
 A_n=\left\{y:x_d-nr<y_d\leq x_d-(n-1)r\right\}\cup
\]
\[ 
\left\{y:x_d+(n-1)r\leq y_d< x_d+nr\right\}.
\]
Then we obtain
\[
 \int_{|y-x|>3r}\frac{|u(y)|^p}{|x-y|^d}dy=\int_{A_0}\frac{|u(y)|^p}{|x-y|^d}dy+\sum^{\infty}_{n=5}\int_{A_n}\frac{|u(y)|^p}{|x-y|^d}dy
\]
\[
 \leq(3r)^{-d}8r\left\|u\right\|^p_{h^p}+2r^{1-d}\left\|u\right\|^p_{h^p}\sum^{\infty}_{n=5}\frac{1}{(n-1)^d}.
\]
Since $d\geq2$, we have $I_2\leq c_1\left\|u\right\|^p_{h^p}$, where $c_1=c_1(\alpha,d,r)<\infty$. Hence 
\[
 \int_{B(x_0,r)}f_1(x)dx\leq\int_{B(x_0,r)}f_2(x)dx\leq (8r+c_1)\left\|u\right\|^p_{h^p}.
\]
Finally we have
\[
 \int_{B(x_0,r/2)}f_1(x)dx=\C_1\int_{B(x_0,r/2)}\int_{|x_0-y|>r}\left(\frac{r^2-|x_0-x|^2}{|x_0-y|^2-r^2}\right)^{\alpha/2}\frac{|u(y)|^p}{|x-y|^d}dydx
\]
\[
\geq\C_1\int_{B(x_0,r/2)}dx\int_{|x_0-y|>r}\left(\frac{r^2/2}{|x_0-y|^2-r^2}\right)^{\alpha/2}\frac{|u(y)|^p}{2^d|x_0-y|^d}dy
\]
\[
=c_2\int_{|x_0-y|>r}P_{B(x_0,r)}(x_0,y)|u(y)|^pdy\geq c_2|u(x_0)|^p,
\]
where $c_2=c_2(\alpha,d,r)$ and the last estimate follows from Jensen's inequality. Therefore, 
$|u(x_0)|\leq [(8r+c_1)/c_2]^{1/p}\left\|u\right\|_{h^p}$. This gives the conclusion of the lemma.
\end{proof}

\begin{remark}\label{rem:Hpbound}
Lemma~\ref{lem:Hpbound} can also be proved using a modified Poisson kernel for the ball, 
described in \cite[p.65]{BogB}. Here we present a different method, which may be of independent interest.
\end{remark}

\noindent The next property of the Poisson kernels $P_{\HH_{\varepsilon}}(x,y)$ is a counterpart of Lemma~\ref{lem:rotation} from the previous section. 
In the present context, the result is a consequence of the translation invariance of the symmetric stable processes.

\begin{lemma}\label{lem:transl}
Fix $\varepsilon>0$. For all $t,s\in\RR$ such that $|s|<\varepsilon<|t|$ and all $\y\in\RR^{d-1}$ we have 
\[
 \int_{\RR^{d-1}}P_{\HH_{\varepsilon}}\left((\x,t),(\y,s)\right)d\x=\int_{\RR^{d-1}}P_{\HH_{\varepsilon}}\left((\y,t),(\x,s)\right)d\x.
\]
Furthermore, the integrals are constant with respect to $\y$.
\end{lemma}

\begin{proof}
We proceed as in the proof of Lemma~\ref{lem:rotation}, using the translation invariance of $X_t$ instead of the rotation invariance.
\end{proof}

\noindent We will now apply the last three lemmas to obtain the relation between $\|u\|_{h^p}$ and $\f_{\HH}[|u|^p]$.

\begin{lemma}\label{lem:harmaj2}
If $u\in h^p_{\alpha}(\HH)$ for a given $p\in[1,\infty)$, then $\f_{\HH}[|u|^p]$ is finite. Furthermore, 
$\f_{\HH}[|u|^p]\in h^1_{\alpha}(\HH)$ and $\|\f_{\HH}[|u|^p]\|_{h^1}=\|u\|^p_{h^p}$.
\end{lemma}

\begin{proof}
Fix $p\in[1,\infty)$ and let $u\in h^p_{\alpha}(\HH)$. Set
\[
 \f_{\varepsilon}[|u|^p](x):=\EE^x\left|u\left(X(\tau_{\HH_{\varepsilon}})\right)\right|^p.
\]
By Fubini theorem we have
\[
 \|\f_{\varepsilon}[|u|^p]_t\|_1=
\int_{\RR^{d-1}}\int_{\left(\overline{\HH}_{\varepsilon}\right)^c}P_{\HH_{\varepsilon}}\left((\x,t),y\right)|u(y)|^pdyd\x
\]
\[
 =\int^{\varepsilon}_{-\varepsilon}\int_{\RR^{d-1}}|u(\y,s)|^p
 \int_{\RR^{d-1}}P_{\HH_{\varepsilon}}\left((\x,t),(\y,s)\right)d\x d\y ds.
\]
By Lemma~\ref{lem:transl} the last term above is equal to
\[
\int^{\varepsilon}_{-\varepsilon}f_{\varepsilon}(t,s)\int_{\RR^{d-1}}|u(\y,s)|^pd\y ds
\]
\[
\leq\left\|u\right\|^p_{h^p}\int^{\varepsilon}_{-\varepsilon}\int_{\RR^{d-1}}P_{\HH_{\varepsilon}}\left((\overline{w},t),(\x,s)\right)d\x ds
=\left\|u\right\|^p_{h^p}.
\]
In the last integral above, $\overline{w}\in\RR^{d-1}$ is arbitrary. Therefore, $\f_{\varepsilon}[|u|^p]$ is finite a.e. on $\left\{x\in\RR^d:x_d=t\right\}$ 
with respect to $(d-1)$-dimensional Lebesgue measure. In view of the Harnack inequality (see \cite[Theorem 1]{BogKK}), 
$\f_{\varepsilon}[|u|^p]$ is finite everywhere on $\HH_{\varepsilon}$, and thus $\alpha$-harmonic on $\HH_{\varepsilon}$. 
By the Jensen inequality, $\f_{\varepsilon}[u]^p\leq\f_{\varepsilon}[|u|^p]$, so $\f_{\varepsilon}[u]$ also 
is finite. Hence, by Lemma~\ref{lem:Hpbound} and Lemma~\ref{lem:meanvalue}, $u$ satisfies the mean value property (\ref{det19})
for every $\varepsilon>0$ and $x\in\HH_{\varepsilon}$, and thus $|u|^p\leq\f_{\varepsilon}[|u|^p]$. Furthermore, if 
$\varepsilon_1\geq\varepsilon_2>0$, then by the strong Markov property,
\[
 \f_{\varepsilon_1}[|u|^p](x)\leq\f_{\varepsilon_1}[\f_{\varepsilon_2}[|u|^p]](x)=\f_{\varepsilon_2}[|u|^p](x).
\]
Hence the limit $v=\lim_{\varepsilon\to0}\f_{\varepsilon}[|u|^p]$ exists. Since 
$\|\f_{\varepsilon}[|u|^p]_t\|_1\leq\left\|u\right\|^p_{h^p}$ for all $t\neq0$ and $\varepsilon>0$, 
the monotone convergence theorem and the Harnack inequality imply that $v$ is $\alpha$-harmonic on $\HH$ and 
$\|v\|_{h^1}\leq\left\|u\right\|^p_{h^p}$. Obviously, $|u|^p\leq v$ and hence $\|v\|_{h^1}=\left\|u\right\|^p_{h^p}$. 
Furthermore, $\f_{\HH}[|u|^p]\leq v$ by Lemma~\ref{lem:majorant1}. Therefore $\f_{\HH}[|u|^p]\in h^1_{\alpha}$ and 
as in the case of $u$ we conclude that $\f_{\HH}[|u|^p]$ satisfies the mean value property (\ref{det19}) for every 
$\varepsilon>0$ and $x\in\HH_{\varepsilon}$. Hence $\f_{\varepsilon}[|u|^p]\leq\f_{\HH}[|u|^p]$ for every 
$\varepsilon>0$, which gives $v\leq\f_{\HH}[|u|^p]$, so the lemma is proved. 
\end{proof}

\noindent From Lemma~\ref{lem:harmaj2} it follows that $h^p_{\alpha}(\HH)\subset\h^p_{\alpha}(\HH)$ for every $p\in[1,\infty)$. 
The next theorem shows that the equality is no longer true in opposite to the previous section.

\begin{theorem}\label{th:hyperplhp}
Let $u$ be $\alpha$-harmonic on $\HH$.
\begin{enumerate}
	\item[\it 1.] $u\in h^1_{\alpha}(\HH)$ if and only if $u=\p_{\HH}[\mu]$ for some $\mu\in\M(\RR^{d-1})$.
	Furthermore, $\mu$ is unique and $\|\p_{\HH}[\mu]\|_{h^1}=\|\mu\|$.

	\item[\it 2.] $u\in h^p_{\alpha}(\HH)$ for a given $p\in(1,\infty]$ if and only if $u=\p_{\HH}[f]$ for some $f\in L^p(\RR^{d-1})$.
	Furthermore, $f$ is unique and $\|\p_{\HH}[f]\|_{h^p}=\|f\|_p$.
\end{enumerate}
\end{theorem}

\begin{proof}
If $u=\p_{\HH}[\mu]$ for some $\mu\in\M(\RR^{d-1})$, then by Lemma~\ref{lem:poissinteg2} (i) we have $\|u\|_{h^1}\leq\|\mu\|$, 
what gives  $u\in h^1_{\alpha}(\HH)$. Moreover, by Lemma~\ref{lem:poissinteg2} (v), $u_t\to\mu$ weak$^*$ as $t\to0$, so
$\|\mu\|\leq \liminf_{t\to0}\|u_t\|_1$. Hence
\[
\|u\|_{h^1}=\lim_{t\to0}\|u_t\|_1=\|\mu\|.
\]
Conversely, suppose that $u\in h^1_{\alpha}(\HH)$. Then by Lemma~\ref{lem:harmaj2}, $\f_{\HH}[u]$ is finite and 
$\|\f_{\HH}[u]\|_{h^1}=\|u\|_{h^1}$. By Theorem~\ref{th:widerhp}, $u(x)=M_{\HH}[\mu](x)+c|x_d|^{\alpha-1}$ 
for a unique pair $(\mu,c)\in\M(\RR^{d-1})\times\RR$ and by Lemma~\ref{lem:majorantform}, 
$\f_{\HH}[u](x)=M_{\HH}[|\mu|](x)+|c||x_d|^{\alpha-1}$. As $|x_d|^{\alpha-1}\notin h^p_{\alpha}(\HH)$ 
for any $p\in[1,\infty]$, we have $c=0$. By (\ref{det12}) we have $u=\p_{\HH}[\nu]$, where $\nu(d\x)=\mu(d\x)/\p_{\HH}(e_d,\x)$. 
Hence $\f_{\HH}[u]=\p_{\HH}[|\nu|]$ and for $t\in\RR\setminus\left\{0\right\}$ the identity (\ref{det7}) gives
\[
 \|\f_{\HH}[u]_t\|_1=\int_{\RR^{d-1}}\int_{\RR^{d-1}}\p_{\HH}((\x,t),\y)|\nu|(d\y)d\x
\]
\[
=\int_{\RR^{d-1}}\int_{\RR^{d-1}}\p_{\HH}((\y,t),\x)d\x|\nu|(d\y)=\|\nu\|,
\]
so $\nu\in\M(\RR^{d-1})$. This gives the first part.

Let now $1<p\leq\infty$. If $f\in L^p(\RR^{d-1})$ and $u=\p_{\HH}[f]$, then from Lemma~\ref{lem:poissinteg2} (ii), (iv) and (vi) 
it follows, in the same way as for $p=1$, that $u\in h^p_{\alpha}(\HH)$ and $\|u\|_{h^p}=\|f\|_p$. 
Conversely, suppose that $u\in h^p_{\alpha}(\HH)$ for a given $p\in(1,\infty)$. Then by Lemma~\ref{lem:harmaj2}, 
$\f_{\HH}[|u|^p]$ is finite and $\|\f_{\HH}[|u|^p]\|_{h^1}=\|u\|^p_{h^p}$. By Theorem~\ref{th:widerhp}, 
$u=\p_{\HH}[f]$ for a unique function $f\in L^p(\RR^{d-1},\omega_{\alpha})$. Furthermore, by Lemma~\ref{lem:majorantform}, 
$\f_{\HH}[|u|^p]=\p_{\HH}[|f|^p]$ and exactly as in the case $p=1$ we obtain $\|\f_{\HH}[|u|^p]_t\|_1=\|f\|^p_p$ for any $t\neq0$,
so $f\in L^p(\RR^{d-1})$. Finally, let $u\in h^{\infty}_{\alpha}(\HH)$. Then $u$ is bounded on $\HH$, and hence $\f_{\HH}[|u|^q]$ 
is finite for every $q\in[1,\infty)$. Therefore, by Theorem~\ref{th:widerhp}, $u=\p_{\HH}[g]$ for a unique measurable function $g$ such that 
$\|g\|_{q,\alpha}=(\f_{\HH}[|u|^q](e_d))^{1/q}\leq\|u\|_{h^{\infty}}$. Hence $\|g\|_{\infty}<\infty$, so the proof is complete.
\end{proof}

\noindent We will now discuss the corresponding version of the Fatou theorem for the $\alpha$-harmonic functions on $\HH$. For $y\in\LL$ and $\beta>0$, the cone $\Gamma_{\beta}(y)$ on $\HH$ is defined as
\[
\Gamma_{\beta}(y)=\left\{x\in\HH:|x-y|<(1+\beta)|x_d|\right\}.
\]
Analogously to the previous section we define the nontangential limits of functions on $\HH$. 
For a function $u$ on $\RR^d$ we also define the Kelvin transform of $u$ by
\[
K_{\alpha}[u](x)=|x|^{\alpha-d}u(Tx),\quad x\neq0,
\]
and the modified Kelvin transform of $u$ by
\[
\widetilde{K}_{\alpha}[u](x)=2^{(d-\alpha)/2}|x+e_d|^{\alpha-d}u(\widetilde{T}x),\quad x\neq-e_d,
\]
where $T,\widetilde{T}$ are the inversions defined in (\ref{det26}) and (\ref{det14}). Let $D\subset\RR^d\setminus\left\{0\right\}$ be open. In view of \cite[Lemma 7]{BogZ}, a function $u$ is $\alpha$-harmonic on $D$ if and only if $K_{\alpha}[u]$ is $\alpha$-harmonic on $TD$. The scaling property and the translation invariance of symmetric stable processes imply that a function $u$ is $\alpha$-harmonic on $\HH'=\HH\setminus\left\{-e_d\right\}$ if and only if $\widetilde{K}_{\alpha}[u]$ is $\alpha$-harmonic on $\widetilde{T}(\HH')=\DD$. In particular, if $u$ is $\alpha$-harmonic on $\HH$, then $\widetilde{K}_{\alpha}[u]$ is $\alpha$-harmonic on $\DD$. Furthermore, since $\widetilde{T}$ is conformal (see \cite[Proposition 7.18]{ABR}), $u$ has a nontangential limit at $y\in\LL$ if and only if $\widetilde{K}_{\alpha}[u]$ has a nontangential limit at $\widetilde{T}y\in\s$. 

\begin{theorem}\label{th:Fatou2}
Suppose $\mu\in\M(\RR^{d-1})$ and let $\mu(d\x)=f(\x)\omega_{\alpha}(d\x)+\nu(d\x)$ be the Lebesgue decomposition of $\mu$ with respect to $\omega_{\alpha}$.
Then $M_{\HH}[\mu]$ has the nontangential limit $f(\y)$ at almost every $y=(\y,0)\in\LL$ with respect to the $(d-1)$-dimensional Lebesgue measure.
\end{theorem}

\begin{proof}
Using the modified Kelvin transform, Theorem~\ref{th:Fatou} and Lemma~\ref{lem:strongconv} we proceed as in \cite[Theorems 7.28 and 7.29]{ABR}.
\end{proof}

\begin{corollary}\label{cor:Fatou2}
Suppose $u$ is $\alpha$-harmonic and nonnegative on $\HH$. Then $u$ has a nontangential limit at almost every $y\in\LL$ 
with respect to the $(d-1)$-dimensional Lebesgue measure.
\end{corollary}

\begin{proof}
The corollary follows from Proposition~\ref{prop:Martin2} and Theorem~\ref{th:Fatou2}.
\end{proof}

\noindent We will finish this section with few examples of $\alpha$-harmonic functions that do not belong to any of the considered Hardy spaces. 
Let $d=2$ and let $u(x)=u(x_1,x_2)=x_1$. Then one can easily check that $\Delta^{\alpha/2}u\equiv0$ if and only if $\alpha\in(1,2)$ 
(for $\alpha\in(0,1]$ the integrals in (\ref{det27}) are not absolutely convergent). 
First obvious observation is that $u\notin h^p_{\alpha}(\HH)$ for any $p\in[1,\infty]$. 
Furthermore, $u\notin h^1_{\alpha}(\DD)$ (and hence $u\notin h^p_{\alpha}(\DD)$ for any $p\in[1,\infty]$) since $u$ is unbounded at 
infinity. Finally, $\int_{\RR}|x_1|\p_{\HH}(e_2,(x_1,0))dx_1=\infty$, so in view of Lemma~\ref{lem:Martinweak} and 
Theorem~\ref{th:widerhp}, $u\notin\h^1_{\alpha}(\HH)$, and therefore $u\notin\h^p_{\alpha}(\HH)$ for any $p\geq1$. To obtain an 
example of a function bounded at infinity we take $K_{\alpha}[u]=x_1|x|^{\alpha-4}$. Since $T\HH=\HH$, $K_{\alpha}[u]$ is 
$\alpha$-harmonic on $\HH$ by \cite[Lemma 7]{BogZ}. Furthermore, $K_{\alpha}$ is linear, preserves nonnegative $\alpha$-harmonic 
functions and $K_{\alpha}[K_{\alpha}[v]]=v$ for any function $v$ on $\HH$. Let $\h_{\alpha}(\HH)$ denote the set of nonnegative 
$\alpha$-harmonic functions on $\HH$. By Theorem~\ref{th:widerhp} and Proposition~\ref{prop:Martin2} we have 
$\h^1_{\alpha}(\HH)=\h_{\alpha}(\HH)-\h_{\alpha}(\HH)$, and hence  $K_{\alpha}[\h^1_{\alpha}(\HH)]=\h^1_{\alpha}(\HH)$. 
Therefore $K_{\alpha}[u]\notin\h^1_{\alpha}(\HH)$. The last example we give is 
$\widetilde{K}_{\alpha}[u](x)=2^{(4-\alpha)/2}x_1|x+e_2|^{\alpha-4}$, 
the modified Kelvin transform of $u$. $\widetilde{K}_{\alpha}[u]$ is $\alpha$-harmonic on $\DD$ and bounded at infinity. 
Furthermore, for $|x+e_2|<1$ we have 
\[
\frac{|x_1|}{|x+e_2|^{4-\alpha}}\geq\frac{|x_1|}{|x+e_2|^{2}}=\left|\Im\left(\frac{1}{1+z}\right)\right|,
\]
where $z=x_2+ix_1$ and $\Im(1/(1+z))$, the imaginary part of $1/(1+z)$, is a well known example of harmonic function 
of the Laplacian not being in the classical Hardy space $h^1$ on the unit ball in $\RR^2$ (see, e.g., \cite[p. 178]{MR1}). 
Since the norm of $h^1_{\alpha}(\DD)$ is analogous to the classical one, we have $\widetilde{K}_{\alpha}[u]\notin h^1_{\alpha}(\DD)$.

\section{Hitting probabilities for relativistic stable processes}\label{sec5}

In this section we will discuss the behavior of a discontinuous process different than the symmetric $\alpha$-stable one 
in context of the hitting probabilities of $\s$ and $\LL$. Namely, we will consider the so-called 
relativistic $\alpha$-stable process. For a given $m>0$ and $\alpha\in(0,2)$, the relativistic $\alpha$-stable 
process $X^m_t$ is defined as a process with the following characteristic function
\[
\EE^xe^{i\xi\cdot(X^m_t-x)}=e^{mt}e^{-t(|\xi|^2+m^{2/\alpha})^{\alpha/2}},\quad x,\xi\in\RR^d,\quad t\geq0.	
\]
$X^m_t$ is a L\'evy process with the infinitesimal generator equal to $m-(-\Delta+m^{2/\alpha})^{\alpha/2}$. 
The potential theory of this process has been widely studied in recent years, see \cite{R, CS2, KS, GR1, KL, GR2, KSV, ByMR1, ByMR2, CKS}. 
The boundary behavior of functions harmonic for $X^m_t$ reminds the one of $\alpha$-harmonic functions on bounded regular domains, 
see \cite{R, CS2, Ki, GR1, KL}. However, in the case of unbounded sets the behavior is slightly different 
(see \cite{GR2, ByMR1, CKS}), and many properties are still unknown, such as, for example, the Martin representation for nonnegative 
functions harmonic for $X^m_t$ (see \cite{KL}, where the result was proved for bounded $\kappa$-fat sets). In this section we aim to 
show that $\s$ and $\LL$ are non-polar for $X^m_t$ for $\alpha\in(1,2)$, what may serve as a motivation for studying 
the boundary value problems for $m-(-\Delta+m^{2/\alpha})^{\alpha/2}$ on $\DD$ and $\HH$.

For a set $B\subset\RR^d$ denote 
$\tau^m_B=\inf\left\{t>0:X^m_t\notin B\right\}$, $T^m_B=\tau^m_{B^c}$. 
Let $S_t$ be the usual $\alpha/2$-stable subordinator given by the Laplace exponent $\psi(\lambda)=\lambda^{\alpha/2}$, and 
let $S^m_t$ be the relativistic $\alpha/2$-stable subordinator, i.e., the subordinator with the Laplace 
exponent $\psi_m(\lambda)=(\lambda+m^{2/\alpha})^{\alpha/2}-m$. Denote by $h(t,x)$, $h_m(t,x)$ the transition densities of
$S_t$ and $S^m_t$, respectively. We have
\begin{equation}\label{det44}
h_m(t,x)=e^{mt}e^{-m^{2/\alpha}x}h(t,x),
\end{equation}
see \cite{R}. It is a well known fact that the potential density of $S_t$ is equal to $x^{\alpha/2-1}/\Gamma(\alpha/2)$, see,
e.g, \cite[p. 97]{BBKRSV}. Furthermore, the potential density of $S^m_t$ is given by
\begin{equation}\label{det39}
q_m(x)=e^{-m^{2/\alpha}x}x^{\alpha/2-1}E_{\alpha/2,\alpha/2}(mx^{\alpha/2}),\quad x>0,
\end{equation}
where 
\[
E_{\gamma,\beta}(t)=\sum^{\infty}_{n=0}\frac{t^n}{\Gamma(\beta+\gamma n)},\quad \gamma,\beta,t>0,
\]
is the two parameter Mittag-Leffler function, see \cite[p. 97]{BBKRSV}. 
Let $Z^m_t:=|X^m_t|$ and let $B_t$ be the usual Brownian motion on $\RR^d$. Set 
$Y(t):=|B_t|$. Then $Z^m_t$ and $Y(S^m_t)$ are equivalent provided that $Y(t)$ and $S^m_t$ are independent. 
It is well known (see, e.g., \cite[p. 116]{P2}) that the transition function of $Y(t)$ is given by 
\[
\PP^x(Y(t)\in A)=\int_A f(t,x,y)\mu(dy),
\] 
where $\mu(dy)=2^{-d/2}[\Gamma(d/2+1)]^{-1}y^ddy$, and 
\begin{equation}\label{det40}
f(t,x,y)=\Gamma(d/2)\frac{1}{2t}\left(\frac{xy}{2}\right)^{1-d/2}\exp\left[\frac{-(x^2+y^2)}{4t}\right]I_{d/2-1}\left(\frac{xy}{2t}\right),
\end{equation}
where $I_{\vartheta}$ is the modified Bessel function of the first kind. This then implies, that the transition function of $Z^m_t$ is given 
by
\[
\PP^x(Z^m_t\in A)=\int_A f_m(t,x,y)\mu(dy),
\] 
where $\mu$ was defined above and
\begin{equation}\label{det41}
f_m(t,x,y)=\int^{\infty}_0h_m(t,s)f(s,x,y)ds.
\end{equation}
Furthermore, the potential density of $Z^m_t$ is given by
\begin{equation}\label{det42}
u_m(x,y)=\int^{\infty}_0q_m(t)f(t,x,y)dt.
\end{equation}
For $r>0$ let $\s_r:=r\s$. It is clear that $\s_r$ is polar for $X^m_t$ if and only if $\left\{r\right\}$ is polar for $Z^m_t$.
Let $\Phi^m_r(x):=\PP^x(T^m_{\s_r}<\infty)$ be the hitting probability of $\s_r$ for $X^m_t$. In the next proposition we consider only 
$d\geq2$.

\begin{proposition}\label{prop:sphererelat}
$\left\{r\right\}$ is polar for $Z^m_t$ provided $\alpha\in(0,1]$. When $\alpha\in(1,2)$ then $r$ is regular for $\left\{r\right\}$ for $Z^m_t$ and 
for all $\tilde{x}\in\RR^d$ we have
\begin{equation}\label{det43}
\Phi^m_r(\tilde{x})=\begin{cases}
		1, & \ d=2,\\
		\frac{u_m(|\tilde{x}|,r)}{u_m(r,r)}, & \ d\geq3.
		\end{cases}
\end{equation}
\end{proposition}

\begin{proof} 
We apply the methods of \cite[proof of Proposition 2.1]{P2}. For $\lambda>0$ we define the $\lambda$-potential of $Z^m_t$ by
\[
u^{\lambda}_m(x,y)=\int^{\infty}_0e^{-\lambda t}f_m(t,x,y)dt,\quad x,y\geq0.
\]
By (\ref{det44}), (\ref{det41}) and Fubini theorem we have
\[
u^{\lambda}_m(x,y)=\int^{\infty}_0e^{-m^{2/\alpha}s}f(s,x,y)\int^{\infty}_0e^{(m-\lambda)t}h(t,s)dtds.
\]
In view of (\ref{det39}), for $0<\lambda<m$ we obtain
\begin{equation}\label{det47}
u^{\lambda}_m(x,y)=\int^{\infty}_0e^{-m^{2/\alpha}s}s^{\alpha/2-1}f(s,x,y)E_{\alpha/2,\alpha/2}((m-\lambda)s^{\alpha/2})ds.
\end{equation}
We have
\[
E_{\gamma,\beta}(t)\cong \gamma^{-1}t^{(2-2\beta)/(2\gamma)}e^{t^{1/\gamma}},\quad t\to\infty,
\]
see, e.g., \cite[Theorem 1]{G}. Here $f(t)\cong g(t)$ means that the ratio of $f$ and $g$ tends to $1$. 
Furthermore, we have the following asymptotics of the Bessel function $I_{\vartheta}$:
\[
I_{\vartheta}(r)\cong\frac{1}{\Gamma(\vartheta+1)}\left(\frac{r}{2}\right)^{\vartheta},\quad r\to0^+,
\]
\[
I_{\vartheta}(r)\cong(2\pi r)^{-1/2}e^r,\quad r\to\infty.
\]
Hence, for $x,y>0$ we get
\[
e^{-m^{2/\alpha}s}f(s,x,y)s^{\alpha/2-1}E_{\alpha/2,\alpha/2}((m-\lambda)s^{\alpha/2})
\]
\begin{equation}\label{det45}
\cong\frac{\Gamma(d/2)}{\Gamma(\alpha/2)}2^{d/2-2}\pi^{-1/2}(xy)^{(1-d)/2}e^{-m^{2/\alpha}s}e^{-(x-y)^2/(4s)}s^{(\alpha-3)/2},\quad s\to0^+.
\end{equation}
When $x=0$ or $y=0$, then
\[
e^{-m^{2/\alpha}s}f(s,x,y)s^{\alpha/2-1}E_{\alpha/2,\alpha/2}((m-\lambda)s^{\alpha/2})
\]
\begin{equation}\label{det49}
\cong\Gamma(\alpha/2)^{-1}2^{-d/2}e^{-m^{2/\alpha}s}e^{-(x^2+y^2)/(4s)}s^{\alpha/2-d/2-1},\quad s\to0^+,
\end{equation}
and for all $x,y\geq0$ we get
\[
e^{-m^{2/\alpha}s}f(s,x,y)s^{\alpha/2-1}E_{\alpha/2,\alpha/2}((m-\lambda)s^{\alpha/2})
\]
\begin{equation}\label{det46}
\cong\alpha^{-1}2^{1-d/2}(m-\lambda)^{2/\alpha-1}e^{[(m-\lambda)^{2/\alpha}-m^{2/\alpha}]s}e^{-(x^2+y^2)/(4s)}s^{-d/2},\quad s\to\infty.
\end{equation}
Therefore, for $\alpha\in(0,1]$ and $y>0$ we have $u^{\lambda}_m(x,y)\to\infty$ 
when $x\to y$, while for $\alpha\in(1,2)$, $u^{\lambda}_m(x,y)$ is bounded and continuous in $x$ in a neighborhood of $y$. For $r>0$ let 
$T^m_r:=\inf\left\{t>0:Z^m_t=r\right\}$. As in the proof of \cite[Proposition 2.1]{P2} we conclude that if $\alpha\in(0,1]$ then $\left\{r\right\}$ is polar 
for $Z^m_t$ for every $r$, while for $\alpha\in(1,2)$ we have
\begin{equation}\label{det48}
\EE^x(e^{-\lambda T^m_r};T^m_r<\infty)=\frac{u^{\lambda}_m(x,r)}{u^{\lambda}_m(r,r)}>0,\quad x\geq0,r>0.
\end{equation}
This then implies that $r$ is regular for $\left\{r\right\}$ for $Z^m_t$. Suppose first $d\geq3$. In view of (\ref{det42}), (\ref{det47}), 
(\ref{det45}), (\ref{det49}) and (\ref{det46}), for all $x\geq0$, $r>0$ and $\alpha\in(1,2)$ we have $u^{\lambda}_m(x,r)\to u_m(x,r)<\infty$ as $\lambda\to0$.
Hence
\[
\Phi^m_r(\tilde{x})=\PP^{|\tilde{x}|}(T^m_r<\infty)=\frac{u_m(|\tilde{x}|,r)}{u_m(r,r)},\quad\tilde{x}\in\RR^d.
\]
Let now $d=2$. Then by (\ref{det46}), $u^{\lambda}_m(x,r)\to u_m(x,r)=\infty$ as $\lambda\to0$, while 
(\ref{det45}), (\ref{det49}) and (\ref{det46}) imply that $u^{\lambda}_m(x,r)/u^{\lambda}_m(r,r)\to1$ as $\lambda\to0$. By (\ref{det48}) 
we get $\Phi^m_r(\tilde{x})=1$ for all $\tilde{x}\in\RR^d$, as desired.
\end{proof}

\begin{remark}\label{rem:hyperrelat}
For $d=1$ the relativistic $\alpha$-stable process is pointwise recurrent (i.e. hits points almost surely) for $\alpha\in(1,2)$ 
and transient for $\alpha\in(0,1]$. This can be proved using similar methods as in the proof of Proposition~\ref{prop:sphererelat}, see also 
\cite{P1} for the symmetric $\alpha$-stable case. Furthermore, from the subordination it follows that the coordinates of $X^m_t$ are one-dimensional 
relativistic $\alpha$-stable motions. As a consequence, for $d\geq2$ and $\alpha\in(1,2)$ the process $X^m_t$ hits $\LL$ almost surely (see formula 
(\ref{det32})).
\end{remark}

We would like to point out that in view of (\ref{det3}), for $d=2$ and $\alpha\in(1,2)$ the hitting probability of 
$\s$ for the symmetric 
$\alpha$-stable process is strictly less than 1, in opposite to the present case. It seems to be an interesting  
problem to find some estimates of the harmonic measures of $\DD$ and $\HH$ for $X^m_t$, $\alpha\in(1,2)$. 
We will give here the explicit formula only for the so-called $\lambda$-harmonic measure of $\HH$ for $X^m_t$ in the 
particular case $\lambda=m$. 

For an open set $D\subset\RR^d$, the $\lambda$-harmonic measure of $D$ for $X^m_t$ is given by
\[
\omega^x_{m,\lambda}(A,D)
:=\EE^x\left[\exp(-\lambda\tau^{m}_D){\bf 1}_A(X^{m}(\tau^{m}_D));\tau^{m}_D<\infty\right],\quad x\in D,A\subset D^c,
\] 
where ${\bf 1}_A$ is the indicator of $A$. For simplicity 
denote $\omega^x_m(A,D):=\omega^x_{m,m}(A,D)$. The measure $\omega^x_{m}$ can be regarded as the harmonic measure for the 
relativistic $\alpha$-stable process killed at an independent exponential time with expectation $1/m$. The infinitesimal 
generator of this process is equal to $-(m^{2/\alpha}-\Delta)^{\alpha/2}$.

\begin{proposition}\label{prop:hyperplrelat}
Let $\alpha\in(1,2)$. The measure $\omega^x_m(\cdot,\HH)$ has a density with respect to the $(d-1)$-dimensional Lebesgue measure on $\LL$ given by
\[
\p^m_{\HH}(x,y)=\C_4\frac{|x_d|^{\alpha-1}}{|x-y|^{\frac{d+\alpha-2}{2}}}K_{\frac{d+\alpha-2}{2}}\left(m^{1/\alpha}|x-y|\right),\quad x\in\HH,y\in\LL,
\]
where $\C_4=((\alpha-1)(m^{1/\alpha}/2)^{(d+\alpha-2)/2})/(\pi^{(d-1)/2}\Gamma((\alpha+1)/2))$ and $K_{\vartheta}$ is the modified 
Bessel function of the third kind.
\end{proposition}

\begin{proof} 
The proof is based on the results of \cite{ByMR2} and is a small modification of the proof of Proposition~\ref{prop:poisker}. 
We use the following integral representation of the function 
$K_{\vartheta}(z)$:
\[
K_{\vartheta}(z)=2^{-\vartheta-1}z^{\vartheta}\int^{\infty}_0e^{-t}e^{-z^2/(4t)}t^{-\vartheta-1}dt,
\]
where $\Re(z^2)>0$, $|\arg z|<\frac{\pi}{2}$, $\vartheta\in\RR$, see \cite[Vol. II, 7.11 (23)]{E}.
\end{proof}

{\bf Acknowledgements.} This paper is a part of the author's Ph.D. thesis. He would like to thank Prof. P. Graczyk, his supervisor, 
for his help and guidance in preparing the manuscript. Special thanks go to Prof. K. Bogdan for many helpful suggestions and remarks, 
and to Prof. M. Ryznar for pointing out his results in \cite{ByMR2}.


\begin{thebibliography}{HD}


\bibitem{A} H. Aikawa,
\newblock {\em Boundary Harnack principle and Martin boundary for a uniform domain}, 
\newblock J. Math. Soc. Japan Vol. 53 (2001), No. 1, 119-145.
	      
	      
\bibitem{An} A. Ancona, 
\newblock {\em Principe de Harnack \`a la fronti\`ere et th\'eor\`eme de Fatou pour un 
op\'erateur elliptique dans un domaine lipschitzien}, 
\newblock Ann. Inst. Fourier (Grenoble) 28 (1978), no. 4, 169-213.


\bibitem{ABR} S. Axler, P. Bourdon, W. Ramey,
\newblock {\em Harmonic Function Theory}, 
\newblock Springer-Verlag, New York, Berlin, 1992.


\bibitem{BaY} R.F. Bass, D. You, 
\newblock {\em A Fatou theorem for $\alpha$-harmonic functions},
\newblock Bull. Sci. Math. 127(7), 635-648 (2003).


\bibitem{Bog} K. Bogdan,
\newblock {\em The boundary Harnack principle for the fractional Laplacian},
\newblock Studia Math. 123 (1997), 43-80.


\bibitem{Bog2} K. Bogdan,
\newblock {\em Representation of $\alpha$-harmonic functions in Lipschitz domains},
\newblock Hiroshima Math. J. 29 (1999), 227-243.


\bibitem{BogB} K. Bogdan, T. Byczkowski, 
\newblock {\em Potential theory for the $\alpha$-stable Schr\"odinger operator on bounded Lipschitz domains},
\newblock Studia Math. 133 (1999), 53-92.


\bibitem{BBKRSV} K. Bogdan, T. Byczkowski, T. Kulczycki, M. Ryznar, R. Song, Z. Vondra{\v{c}}ek,
\newblock {\em Potential analysis of stable processes and its extensions},
  volume 1980 of {\em Lecture Notes in Mathematics}.
\newblock Springer-Verlag, Berlin, 2009.


\bibitem{BD} K. Bogdan, B. Dyda,
\newblock {\em Relative Fatou theorem for harmonic functions of rotation invariant stable processes in smooth domains}, 
\newblock Studia Math. 157 (2003), 83-96. MR 1 980 119.


\bibitem{BDL} K. Bogdan, B. Dyda, T. Luks,
\newblock {\em On Hardy spaces},
\newblock preprint 2011, available at http://arxiv.org/abs/1109.0210


\bibitem{BJ}  K. Bogdan, T. Jakubowski, 
\newblock {\em Probleme de Dirichlet pour les fonctions $\alpha$-harmoniques sur les domaines coniques},
\newblock Ann. Math. Blaise Pascal 12 (2005) 297-308.


\bibitem{BogKK}  K. Bogdan, T. Kulczycki, M. Kwa\'snicki, 
\newblock {\em Estimates and structure of $\alpha$-harmonic functions},
\newblock Probab. Theory Relat. Fields 140(3-4), 345-381 (2008).


\bibitem{BogZ}  K. Bogdan, T. \.Zak, 
\newblock {\em On Kelvin transformation},
\newblock J. Theor. Prob. 19(1), 89-120 (2006).


\bibitem{BorS}  A. N. Borodin, P. Salminen, 
\newblock {\em Handbook of Brownian Motion - Facts and Formulae},
\newblock Birkh\"auser Verlag, Basel, 2 edition, 2002.


\bibitem{ByMR1} T. Byczkowski, J. Ma\l{}ecki, M. Ryznar,
\newblock {\em Bessel potentials, hitting distributions and Green functions},
\newblock Trans. Amer. Math. Soc. 361 (2009), no. 9, 4871-4900. 


\bibitem{ByMR2} T. Byczkowski, J. Ma\l{}ecki, M. Ryznar,
\newblock {\em Hitting Half-spaces by Bessel-Brownian Diffusions},
\newblock Pot. Anal. Volume 33, Number 1 (2010), 47-83.


\bibitem{C} L. Caffarelli, 
\newblock {\em Variational problems for free boundaries for the fractional Laplacian},
\newblock J. Eur. Math. Soc. 12 (2010), no. 5, 1151-1179.


\bibitem{CKS} Z.-Q. Chen, P. Kim, R. Song,
\newblock {\em Green function estimates for relativistic stable processes in half-space-like open sets},
\newblock Stochastic Process. Appl. 121 (2011), no. 5, 1148-1172. 


\bibitem{CS1} Z.-Q. Chen, R. Song,
\newblock {\em Martin boundary and integral representation for harmonic functions of symmetric stable processes},
\newblock J. Funct. Anal. 159 (1998), no. 1, 267-294. 


\bibitem{CS2} Z.-Q. Chen, R. Song,
\newblock {\em Drift transforms and Green function estimates for discontinuous processes},
\newblock J. Funct. Anal. 201 (2003), no. 1, 262-281.


\bibitem{Do2} J.L. Doob, 
\newblock {\em Classical Potential Theory and Its Probabilistic Counterpart},
\newblock Springer-Verlag, New York, 1984.


\bibitem{D} P.L. Duren,
\newblock {\em Theory of $H^p$ Spaces},
\newblock Academic Press, New York, 1970; reprinted with supplement by Dover Publications, Mineola, N.Y., 2000.


\bibitem{E} A. Erdelyi,
\newblock {\em Higher transcendental functions}, Bateman Manuscript Project, 
\newblock Vol. I,II, McGraw-Hill, New York, 1953.


\bibitem{G} S. Gerhold,
\newblock {\em Asymptotics for a variant of the Mittag-Leffler function},
\newblock preprint 2011, available at http://arxiv.org/abs/1103.2285


\bibitem{GR1} T. Grzywny, M. Ryznar,
\newblock {\em Estimates of Green function for some perturbations of fractional Laplacian},
\newblock Illinois J. Math 51 (2007), no. 4, 1409-1438.


\bibitem{GR2} T. Grzywny, M. Ryznar,
\newblock {\em Two-sided optimal bounds for half-spaces for relativistic $\alpha$-stable process},
\newblock Pot. Anal. 28 (2008), no. 3, 201-239.
              

\bibitem{I} Y. Isozaki,
\newblock {\em Hitting of a line or a half-line in the plane by two-dimensional symmetric stable L\'evy processes}, 
\newblock Stochastic Process. Appl. 121 (2011) 1749-1769.
              
              
\bibitem{JK} D.S. Jerison, C.E. Kenig, 
\newblock {\em Boundary value problems on Lipschitz domains},
\newblock In {\em Studies in partial differential equations}, volume 23 of {\em
  MAA Stud. Math.}, Math. Assoc. America, Washington, DC, 1982, 1-68.

\bibitem{Ki} P. Kim,
\newblock {\em Relative Fatou's theorem for  $(-\Delta)^{\alpha/2}$-harmonic functions in bounded $\kappa$-fat open sets}, 
\newblock J. Funct. Anal. 234 (2006), no. 1, 70-105.


\bibitem{KL} P. Kim, Y.-R. Lee,
\newblock {\em Generalized 3G theorem and application to relativistic stable process on non-smooth open sets}, 
\newblock J. Funct. Anal. 246 (2007), no. 1, 113-143.


\bibitem{KSV} P. Kim, R. Song, Z. Vondra{\v{c}}ek,
\newblock {\em Boundary Harnack principle for subordinate Brownian motions}, 
\newblock Stochastic Process. Appl. 119 (2009), no. 5, 1601-1631. 


\bibitem{K}  P. Koosis,
\newblock {\em Introduction to $H_p$ Spaces},
\newblock Cambridge University Press, Cambridge 1998.


\bibitem{KS}  T. Kulczycki, B. Siudeja,
\newblock {\em Intrinsic ultracontractivity of the Feynman-Kac semigroup for relativistic stable processes},
\newblock Trans. Amer. Math. Soc. 358 (2006), no. 11, 5025-5057. 


\bibitem{L} T. Luks,
\newblock {\em Hardy spaces for the Laplacian with lower order perturbations},
\newblock Studia Math. 204 (2011), 39-62.


\bibitem{MR} K. Michalik, M. Ryznar, 
\newblock {\em Relative Fatou theorem for $\alpha$-harmonic functions in Lipschitz domains}, 
\newblock Illinois J. Math. 48 (2004), no. 3, 977-998.


\bibitem{MR1} K. Michalik, M. Ryznar, 
\newblock {\em Hardy spaces for $\alpha$-harmonic functions in regular domains},
\newblock Math. Z., Volume 265, Number 1 (2010), 173-186.


\bibitem{MS} K. Michalik, K. Samotij, 
\newblock {\em Martin representation for $\alpha$-harmonic functions},
\newblock Probab. Math. Stat. 20 (2000), 75-91.


\bibitem{P1} S.C. Port, 
\newblock {\em Hitting times and potential for recurrent stable processes},
\newblock J. Analyse Math. 20 (1967), 371-395.


\bibitem{P2} S.C. Port, 
\newblock {\em The First Hitting Distribution of a Sphere for Symmetric Stable Processes},
\newblock Trans. Amer. Math. Soc., 135 (1969), 115-125.


\bibitem{R} M. Ryznar, 
\newblock {\em Estimates of Green Function for Relativistic $\alpha$-Stable Process},
\newblock Pot. Anal. Volume 17, Number 1 (2002), 1-23.


\bibitem{SW}  R. Song, J.-M. Wu,
\newblock {\em Boundary Harnack principle for symmetric stable processes},
\newblock J. Funct. Anal. 168 (1999), 403-427.


\bibitem{S1}  E.M. Stein, 
\newblock {\em Singular integrals and differentiability properties of functions},
\newblock Princeton University Press, Princeton 1970.


\bibitem{S2} E.M. Stein, 
\newblock {\em Boundary behavior of holomorphic functions of several complex variables}, 
\newblock Princeton University Press and University of Tokyo Press, Princeton, New Jersey, 1972.


\bibitem{W} K.O. Widman, 
\newblock {\em On the boundary behavior of solutions to a class of elliptic partial differential equations}, 
\newblock Ark. Mat. 6 (1966), 485-533.


\bibitem{Wu1} J.-M. Wu, 
\newblock {\em Comparisons of kernel functions, boundary Harnack principle and relative Fatou theorem on Lipschitz domains},
\newblock Ann. Inst. Fourier Grenoble 28 (1978), 147-167.


\bibitem{Wu2} J.-M. Wu, 
\newblock {\em Harmonic measures for symmetric stable processes}, 
\newblock Studia Math. 149 (2002), 281-293.


\end{thebibliography}
\end{document}